\documentclass[12pt]{article}
\usepackage{full page, 
amssymb, amscd, graphicx, 
}
    \title{{\bf Logarithmic
tensor category theory, II: Logarithmic formal calculus
and properties of logarithmic intertwining operators}}
    \author{Yi-Zhi Huang, James Lepowsky and Lin Zhang}
    \date{}

\begin{document}
    \bibliographystyle{alpha}
    \maketitle

    \newtheorem{rema}{Remark}[section]
    \newtheorem{propo}[rema]{Proposition}
    \newtheorem{theo}[rema]{Theorem}
   \newtheorem{defi}[rema]{Definition}
    \newtheorem{lemma}[rema]{Lemma}
    \newtheorem{corol}[rema]{Corollary}
     \newtheorem{exam}[rema]{Example}
\newtheorem{assum}[rema]{Assumption}
     \newtheorem{nota}[rema]{Notation}
        \newcommand{\ba}{\begin{array}}
        \newcommand{\ea}{\end{array}}
        \newcommand{\be}{\begin{equation}}
        \newcommand{\ee}{\end{equation}}
        \newcommand{\bea}{\begin{eqnarray}}
        \newcommand{\eea}{\end{eqnarray}}
        \newcommand{\nno}{\nonumber}
        \newcommand{\nn}{\nonumber\\}
        \newcommand{\lbar}{\bigg\vert}
        \newcommand{\p}{\partial}
        \newcommand{\dps}{\displaystyle}
        \newcommand{\bra}{\langle}
        \newcommand{\ket}{\rangle}
 \newcommand{\res}{\mbox{\rm Res}}
\newcommand{\wt}{\mbox{\rm wt}\;}
\newcommand{\swt}{\mbox{\scriptsize\rm wt}\;}
 \newcommand{\pf}{{\it Proof}\hspace{2ex}}
 \newcommand{\epf}{\hspace{2em}$\square$}
 \newcommand{\epfv}{\hspace{1em}$\square$\vspace{1em}}
        \newcommand{\ob}{{\rm ob}\,}
        \renewcommand{\hom}{{\rm Hom}}
\newcommand{\C}{\mathbb{C}}
\newcommand{\R}{\mathbb{R}}
\newcommand{\Z}{\mathbb{Z}}
\newcommand{\N}{\mathbb{N}}
\newcommand{\A}{\mathcal{A}}
\newcommand{\Y}{\mathcal{Y}}
\newcommand{\Arg}{\mbox{\rm Arg}\;}
\newcommand{\comp}{\mathrm{COMP}}
\newcommand{\lgr}{\mathrm{LGR}}

\newcommand{\dlt}[3]{#1 ^{-1}\delta \bigg( \frac{#2 #3 }{#1 }\bigg) }

\newcommand{\dlti}[3]{#1 \delta \bigg( \frac{#2 #3 }{#1 ^{-1}}\bigg) }

 \makeatletter
\newlength{\@pxlwd} \newlength{\@rulewd} \newlength{\@pxlht}
\catcode`.=\active \catcode`B=\active \catcode`:=\active \catcode`|=\active
\def\sprite#1(#2,#3)[#4,#5]{
   \edef\@sprbox{\expandafter\@cdr\string#1\@nil @box}
   \expandafter\newsavebox\csname\@sprbox\endcsname
   \edef#1{\expandafter\usebox\csname\@sprbox\endcsname}
   \expandafter\setbox\csname\@sprbox\endcsname =\hbox\bgroup
   \vbox\bgroup
  \catcode`.=\active\catcode`B=\active\catcode`:=\active\catcode`|=\active
      \@pxlwd=#4 \divide\@pxlwd by #3 \@rulewd=\@pxlwd
      \@pxlht=#5 \divide\@pxlht by #2
      \def .{\hskip \@pxlwd \ignorespaces}
      \def B{\@ifnextchar B{\advance\@rulewd by \@pxlwd}{\vrule
         height \@pxlht width \@rulewd depth 0 pt \@rulewd=\@pxlwd}}
      \def :{\hbox\bgroup\vrule height \@pxlht width 0pt depth
0pt\ignorespaces}
      \def |{\vrule height \@pxlht width 0pt depth 0pt\egroup
         \prevdepth= -1000 pt}
   }
\def\endsprite{\egroup\egroup}
\catcode`.=12 \catcode`B=11 \catcode`:=12 \catcode`|=12\relax
\makeatother

\def\hboxtr{\FormOfHboxtr} 
\sprite{\FormOfHboxtr}(25,25)[0.5 em, 1.2 ex] 

:BBBBBBBBBBBBBBBBBBBBBBBBB |
:BB......................B |
:B.B.....................B |
:B..B....................B |
:B...B...................B |
:B....B..................B |
:B.....B.................B |
:B......B................B |
:B.......B...............B |
:B........B..............B |
:B.........B.............B |
:B..........B............B |
:B...........B...........B |
:B............B..........B |
:B.............B.........B |
:B..............B........B |
:B...............B.......B |
:B................B......B |
:B.................B.....B |
:B..................B....B |
:B...................B...B |
:B....................B..B |
:B.....................B.B |
:B......................BB |
:BBBBBBBBBBBBBBBBBBBBBBBBB |

\endsprite

\def\shboxtr{\FormOfShboxtr} 
\sprite{\FormOfShboxtr}(25,25)[0.3 em, 0.72 ex] 

:BBBBBBBBBBBBBBBBBBBBBBBBB |
:BB......................B |
:B.B.....................B |
:B..B....................B |
:B...B...................B |
:B....B..................B |
:B.....B.................B |
:B......B................B |
:B.......B...............B |
:B........B..............B |
:B.........B.............B |
:B..........B............B |
:B...........B...........B |
:B............B..........B |
:B.............B.........B |
:B..............B........B |
:B...............B.......B |
:B................B......B |
:B.................B.....B |
:B..................B....B |
:B...................B...B |
:B....................B..B |
:B.....................B.B |
:B......................BB |
:BBBBBBBBBBBBBBBBBBBBBBBBB |

\endsprite


\begin{abstract}
This is the second part in a series of papers in which 
we introduce and develop a natural, general tensor category theory for
suitable module categories for a vertex (operator) algebra.  
In this paper (Part II), we develop logarithmic formal calculus 
and study logarithmic intertwining
operators.
\end{abstract}

\vspace{2em}

In this paper, Part II of a series of eight papers on logarithmic
tensor category theory, we develop logarithmic formal calculus and
study logarithmic intertwining operators.  The sections, equations,
theorems and so on are numbered globally in the series of papers
rather than within each paper, so that for example equation (a.b) is
the b-th labeled equation in Section a, which is contained in the
paper indicated as follows: In Part I \cite{HLZ1}, which contains
Sections 1 and 2, we give a detailed overview of our theory, state our
main results and introduce the basic objects that we shall study in
this work.  We include a brief discussion of some of the recent
applications of this theory, and also a discussion of some recent
literature.  The present paper, Part II, contains Section 3.  In Part
III \cite{HLZ3}, which contains Section 4, we introduce and study
intertwining maps and tensor product bifunctors.  In Part IV
\cite{HLZ4}, which contains Sections 5 and 6, we give constructions of
the $P(z)$- and $Q(z)$-tensor product bifunctors using what we call
``compatibility conditions'' and certain other conditions.  In Part V
\cite{HLZ5}, which contains Sections 7 and 8, we study products and
iterates of intertwining maps and of logarithmic intertwining
operators and we begin the development of our analytic approach.  In
Part VI \cite{HLZ6}, which contains Sections 9 and 10, we construct
the appropriate natural associativity isomorphisms between triple
tensor product functors.  In Part VII \cite{HLZ7}, which contains
Section 11, we give sufficient conditions for the existence of the
associativity isomorphisms.  In Part VIII \cite{HLZ8}, which contains
Section 12, we construct braided tensor category structure.

\paragraph{Acknowledgments}
The authors gratefully
acknowledge partial support {}from NSF grants DMS-0070800 and
DMS-0401302.  Y.-Z.~H. is also grateful for partial support {}from NSF
grant PHY-0901237 and for the hospitality of Institut des Hautes 
\'{E}tudes Scientifiques in the fall of 2007.

\renewcommand{\theequation}{\thesection.\arabic{equation}}
\renewcommand{\therema}{\thesection.\arabic{rema}}
\setcounter{section}{2}
\setcounter{equation}{0}
\setcounter{rema}{0}

\section{Logarithmic formal calculus, logarithmic intertwining
operators and their basic properties}

In this section we study the notion of ``logarithmic intertwining
operator'' introduced in \cite{Mi}.  For this, we will need to discuss
spaces of formal series in powers of both $x$ and ``$\log x$'', a new
formal variable, with coefficients in certain vector spaces.  In this
``logarithmic formal calculus,'' we establish certain properties, some
of them quite subtle, of the formal derivative operator $d/dx$ acting
on such spaces. Then, following \cite{Mi} with a slight variant (see
Remark \ref{log:compM}), we introduce the notion of logarithmic
intertwining operator. These are the appropriate replacements of
ordinary intertwining operators when $L(0)$-semisimplicity is relaxed.
In the strongly graded setting, it will be natural to consider the
associated ``grading-compatible'' logarithmic intertwining operators.
We work out some important principles and formulas concerning
logarithmic intertwining operators, certain of which turn out to be
the same as in the ordinary intertwining operator case.  Some of these
results require proofs that are quite delicate.

Recall the notation ${\cal W}\{x\}$
(\ref{formalserieswithcomplexpowers}) for the space of formal series
in a formal variable $x$ with coefficients in a vector space ${\cal
W}$, with arbitrary complex powers of $x$.

{}From now on we will sometimes need and use new independent
(commuting) formal variables called $\log x$, $\log y$, $\log x_1$,
$\log x_2, \dots,$ etc. We will work with formal series in such formal
variables together with the ``usual'' formal variables $x$, $y$,
$x_1$, $x_2, \dots,$ etc., with coefficients in certain vector spaces,
and the powers of the monomials in {\it all} the variables can be
arbitrary complex numbers.  (Later we will restrict our attention to
only nonnegative integral powers of the ``$\log$'' variables.)

Given a formal variable $x$, we will use the notation $\frac d{dx}$ to
denote the linear map on ${\cal W}\{x,\log x\}$, for any vector space
${\cal W}$ not involving $x$, defined (and indeed well defined) by the
(expected) formula
\begin{eqnarray}\label{ddxdef}
\lefteqn{\frac d{dx}\bigg(\sum_{m,n\in {\mathbb C}}w_{n,m}x^n(\log
x)^m\bigg)}\nno\\
&&\;\;\;=\sum_{m,n\in{\mathbb C}}((n+1)w_{n+1,m}+ (m+1)w_{n+1,m+1})x^{n}(\log
x)^{m}\nno\\
&&\bigg(=\sum_{m,n\in{\mathbb C}}nw_{n,m}x^{n-1}(\log x)^m+ \sum_{m,n\in
{\mathbb C}}mw_{n,m}x^{n-1}(\log x)^{m-1}\bigg)
\end{eqnarray}
where $w_{n,m}\in {\cal W}$ for all $m, n\in {\mathbb C}$. We will
also use the same notation for the restriction of $\frac d{dx}$ to any
subspace of ${\cal W}\{x,\log x\}$ which is closed under $\frac
d{dx}$, e.g., ${\cal W}\{x\}[[\log x]]$ or
${\mathbb C}[x,x^{-1},\log x]$. Clearly, $\frac
d{dx}$ acting on ${\cal W}\{x\}$ coincides with the usual formal
derivative.

\begin{rema}\label{ddxchk}{\rm
Let $f$, $g$ and $f_i$, $i$ in some index set $I$, all be formal
series of the form
\begin{equation}\label{log:f}
\sum_{m,n\in {\mathbb C}}w_{n,m}x^n(\log x)^m\in {\cal W}\{x,\log x\},\;\;
w_{n,m}\in {\cal W}.
\end{equation}
One checks the following straightforwardly:
Suppose that the sum of $f_i$, $i\in I$, exists (in the obvious
sense). Then the sum of the $\frac d{dx}f_i$,
$i\in I$, also exists and is equal to $\frac d{dx}\sum_{i\in I}f_i$.
More generally, for any $T=p(x)\frac{d}{dx}$, $p(x)\in {\mathbb
C}[x,x^{-1}]$, the sum of $Tf_i$, $i\in I$, exists and is equal to
$T\sum_{i\in I}f_i$. Thus the sum of $e^{yT}f_i$, $i\in I$, exists and is
equal to $e^{yT}\sum_{i\in I}f_i$ ($e^{yT}$ being the
formal exponential series, as usual).  Suppose that ${\cal W}$ is an
(associative) algebra or that the coefficients of either $f$ or $g$
are complex numbers. If the product of $f$ and $g$
exists, then the product of $\frac d{dx}f$ and $g$ and the product of
$f$ and $\frac d{dx}g$ both exist, and $\frac d{dx}(fg)=(\frac
d{dx}f)g+f(\frac d{dx}g)$. Furthermore, for any $T$ as before, the
product of $Tf$ and $g$ and the product of $f$ and $Tg$ both exist,
and $T(fg)=(Tf)g+f(Tg)$. In addition, the product of $e^{yT}f$ and
$e^{yT}g$ exists and is equal to $e^{yT}(fg)$, just as in formulas
(8.2.6)--(8.2.10) of \cite{FLM2}.  The point here, of
course, is just the formal derivation property of $\frac d{dx}$,
except that sums and products of expressions do not exist in general.
}
\end{rema}

\begin{rema}{\rm
Note that the ``equality'' $x=e^{\log x}$ does not hold, since the
left-hand side is a formal variable, while the right-hand side
is a formal series in another formal variable. In fact, this formula
should not be assumed, since, for example,
the formal delta function $\delta(x)=\sum_{n\in {\mathbb Z}}x^n$ would not exist
in the sense of formal calculus, if $x$ were allowed to be replaced by
the formal series $e^{\log x}$. By contrast, note that the equality
\begin{equation}\label{log:logex}
\log e^x=x
\end{equation}
does indeed hold. This is because the formal series $e^x$ is of the
form $1+X$ where $X$ involves only positive integral powers of
$x$ and in (\ref{log:logex}), ``$\log$'' refers to the usual formal
logarithmic series
\begin{equation}\label{log:usual}
\log(1+X)=\sum_{i\geq 1} \frac{(-1)^{i-1}}i X^i,
\end{equation}
{\em not} to the ``$\log$'' of a formal variable. We will use the
symbol ``$\log$'' in both ways, and the meaning will be clear in
context.  }
\end{rema}

We will typically use notations of the form $f(x)$, instead of
$f(x,\log x)$, to denote elements of ${\cal W}\{x,\log x\}$ for some
vector space ${\cal W}$ as above. For this reason, we need to
interpret carefully the meaning of symbols such as $f(x+y)$, or more
generally, symbols obtained by replacing $x$ in $f(x)$ by something
other than just a single formal variable (since $\log x$ is a formal
variable and not the image of some operator acting on $x$). For the
three main types of cases that will be encountered in this work, we
use the following notational conventions; the existence of the
expressions will be justified in Remark \ref{log:exist}:

\begin{nota}
{\rm For formal variables $x$ and $y$, and $f(x)$ of the form
(\ref{log:f}), we define
\begin{eqnarray}\label{log:not1}
f(x+y)&=&\sum_{m,n\in {\mathbb C}}w_{n,m}(x+y)^n\bigg(\log x+
\log\bigg(1+\frac yx\bigg)\bigg)^m\nno\\
&=&\sum_{m,n\in {\mathbb C}}w_{n,m}(x+y)^n\bigg(\log x+
\sum_{i\geq 1} \frac{(-1)^{i-1}}i \bigg(\frac yx\bigg)^i\bigg)^m
\end{eqnarray}
(recall (\ref{log:usual})); in the right-hand side, $(\log x+\sum_{i\geq
1}\frac{(-1)^{i-1}}i
(\frac yx)^i)^m$, according to the binomial expansion convention, is to be
expanded in nonnegative integral powers of the second summand
$\sum_{i\geq 1} \frac{(-1)^{i-1}}i (\frac yx)^i$, so the right-hand
side of (\ref{log:not1}) is equal to
\begin{equation}\label{log:1-tmp}
\sum_{m,n\in {\mathbb C}}w_{n,m}(x+y)^n\sum_{j\in {\mathbb N}}{m\choose j}
(\log x)^{m-j}\bigg(\sum_{i\geq 1} \frac{(-1)^{i-1}}i \bigg(\frac
yx\bigg)^i\bigg)^j
\end{equation}
when expanded one step further. Also define
\begin{equation}\label{log:not2}
f(xe^y)=\sum_{m,n\in {\mathbb C}}w_{n,m}x^ne^{ny}(\log x+y)^m
\end{equation}
and
\begin{equation}\label{log:not3}
f(xy)=\sum_{m,n\in {\mathbb C}}w_{n,m}x^ny^n(\log x+\log y)^m,
\end{equation}
where the binomial expansion convention is again of course being
used.}
\end{nota}

\begin{rema}\label{log:exist}{\rm
The existence of the right-hand side of (\ref{log:not1}), or
(\ref{log:1-tmp}), can be seen by writing $(x+y)^n$ as $x^n(1+\frac
yx)^n$ and observing that
\[
\bigg(\sum_{i\geq 1} \frac{(-1)^{i-1}}i \bigg(\frac yx\bigg)^i\bigg)^j
\in \bigg(\frac yx\bigg)^j{\mathbb C}\bigg[\bigg[\frac yx\bigg]\bigg].
\]
The existence of the right-hand sides of (\ref{log:not2}) and of
(\ref{log:not3}) is clear. Furthermore, both $f(x+y)$ and $f(xe^y)$
lie in ${\cal W}\{x,\log x\}[[y]]$, while $f(xy)$ lies in ${\cal
W}\{xy,\log x\}[[\log y]]$. One might expect that $f(x+y)$ can be
written as $e^{y\frac{d}{dx}}f(x)$, and $f(xe^y)$ as
$e^{yx\frac{d}{dx}}f(x)$ (cf.\ Section 8.3 of \cite{FLM2}), but these
formulas must be verified (see Theorem \ref{log:ids} below).
}
\end{rema}

\begin{rema}\label{subchk}{\rm
It is clear that when there is no $\log x$ involved in $f(x)$, the
expression $f(x+y)$ (respectively, $f(xe^y)$, $f(xy)$) coincides with
the usual formal operation of substitution of $x+y$ (respectively,
$xe^y$, $xy$) for $x$ in $f(x)$. In general, it is straightforward to check
that if the sum of $f_i(x)$, $i\in I$, exists and is equal to $f(x)$,
then the sum of $f_i(x+y)$, $i\in I$ (respectively, $f_i(xe^y), i\in I$,
$f_i(xy), i\in I$), also exists and is equal to $f(x+y)$ (respectively,
$f(xe^y)$, $f(xy)$). Also, suppose that ${\cal W}$ is an (associative)
algebra or that the coefficients of either $f$ or $g$ are complex numbers.
If the product of $f(x)$ and $g(x)$ exists, then
the product of $f(x+y)$ and $g(x+y)$ (respectively, $f(xe^y)$ and
$g(xe^y)$, $f(xy)$ and $g(xy)$) also exists and is equal to
$(fg)(x+y)$ (respectively, $(fg)(xe^y)$, $(fg)(xy)$).  }
\end{rema}

Note that by (\ref{log:not1}),
\begin{equation}\label{logx+y}
\log(x+y)=\log x+\sum_{i\geq 1} \frac{(-1)^{i-1}}i \bigg(\frac
yx\bigg)^i=e^{y\frac d{dx}}\log x.
\end{equation}
The next result includes a generalization of this to arbitrary
elements of ${\cal W}\{x,\log x\}$.  Formula (\ref{log:ck1}) is a
formal ``Taylor theorem'' for logarithmic formal series.  In the case
of non-logarithmic formal series, this principle is used extensively
in vertex operator algebra theory; recall (\ref{formalTaylortheorem})
above and see Proposition 8.3.1 of \cite{FLM2} for the proof in the
generality of formal series with arbitrary complex powers of the
formal variable.  In the logarithmic case below, a much more elaborate
proof is required than in the non-logarithmic case.  The other formula
in Theorem \ref{log:ids}, formula (\ref{log:ck2}), is easier to prove.
It too is important (in the non-logarithmic case) in vertex operator
algebra theory; again see Proposition 8.3.1 of \cite{FLM2}.

\begin{theo}\label{log:ids}
For $f(x)$ as in (\ref{log:f}), we have
\begin{equation}\label{log:ck1}
e^{y\frac d{dx}}f(x)=f(x+y)
\end{equation}
(``Taylor's theorem'' for logarithmic formal series) and
\begin{equation}\label{log:ck2}
e^{yx\frac d{dx}}f(x)=f(xe^y).
\end{equation}
\end{theo}
\pf By Remarks \ref{ddxchk} and \ref{subchk} (or \ref{log:exist}), we
need only prove these equalities for $f(x)=x^n$ and $f(x)=(\log x)^m$,
$m,n\in {\mathbb C}$.  The case $f(x)=x^n$ easily follows {}from the
direct expansion of the two sides of (\ref{log:ck1}) and of
(\ref{log:ck2}) (see Proposition 8.3.1 in \cite{FLM2}). Now assume
that $f(x)=(\log x)^m$, $m\in {\mathbb C}$.

Formula (\ref{log:ck2}) is easier, so we prove it first. By
(\ref{ddxdef}) we have
\[
x\frac d{dx}(\log x)^m=m(\log x)^{m-1},
\]
so that for $k\in {\mathbb N}$,
\[
\bigg(x\frac d{dx}\bigg)^k(\log x)^m=m(m-1)\cdots(m-k+1)(\log
x)^{m-k}=k!{m\choose k}(\log x)^{m-k}.
\]
Thus
\[
e^{yx\frac d{dx}}(\log x)^m=\sum_{k\in {\mathbb N}}
\frac{y^k}{k!}\bigg(x\frac{d}{dx}\bigg)^k(\log x)^m=\sum_{k\in {\mathbb
N}} {m\choose k}y^k (\log x)^{m-k}=(\log x+y)^m,
\]
as we want.

For (\ref{log:ck1}), we shall give two proofs --- an analytic proof and an
algebraic proof. First, consider the analytic function $(\log
z)^m=e^{m\log\log z}$ over, say, $|z-3|<1$ in the complex plane.
In this proof we take the branch of $\log z$ so that
\begin{equation}\label{log:br1}
-\pi< \Im{(\log z)}\leq\pi.
\end{equation}
Then by analyticity, for any $z$ in
this domain, when $|z_1|$ is small enough the Taylor series expansion
$e^{z_1\frac{d}{dz}}(\log z)^m$ converges absolutely to
$(\log(z+z_1))^m$. That is,
\begin{equation}\label{log:ana1}
(\log(z+z_1))^m=e^{z_1\frac{d}{dz}}(\log z)^m
=e^{\frac{z_1}{z}z\frac{d}{dz}}(\log z)^m.
\end{equation}
Observe that as a formal series, the right-hand side of
(\ref{log:ana1}) is in the space $(\log z)^m{\mathbb C}[(\log
z)^{-1}][[z_1/z]]$.

On the other hand, by the choice of domain and the branch of $\log$ we
have
\[
\log(z+z_1)=\log z+\log(1+z_1/z)
\]
and
\[
|\log z|>\log 2>|\log(1+z_1/z)|
\]
when $|z_1|$ is small enough. So when $|z_1|$ is small enough we have
\begin{eqnarray}\label{log:ana2}
(\log(z+z_1))^m&=&(\log z+\log(1+z_1/z))^m\nno\\
&=& \sum_{j\in {\mathbb N}}{m\choose j} (\log z)^{m-j}\bigg(\sum_{i\geq
1} \frac{(-1)^{i-1}}i \bigg(\frac{z_1}z\bigg)^i\bigg)^j.
\end{eqnarray}
Since as formal series, the right-hand sides of (\ref{log:ana1}) and
(\ref{log:ana2}) are both in the space
\[
(\log z)^m{\mathbb C}[(\log z)^{-1}][[z_1/z]],
\]
and both converge to the same analytic function
$(\log(z+z_1))^m$ in the above domain, by setting $z_1=0$ in these two
functions and their derivatives with respect to $z_1$ we see that
their corresponding coefficients of powers of $z_1/z$ and further, of
all monomials in $\log z$ and $z_1/z$ must be the same. Hence we can
replace $z$ and $z_1$ by formal variables $x$ and $y$, respectively,
and obtain (\ref{log:ck1}) for $f(x)=(\log x)^m$.

An algebraic proof of (\ref{log:ck1}) (for $(\log x)^m$) can be given
as follows: Since
\[
\frac d{dx}(\log x)^m=mx^{-1}(\log x)^{m-1}
\]
and higher derivatives involve derivatives of products of powers of
$x$ and powers of $\log x$, let us first compute $(d/dx)^k(x^n(\log
x)^m)$ directly for all $m,n\in {\mathbb C}$ and $k\in{\mathbb N}$. Define
linear maps $T_0$ and $T_1$ on ${\mathbb C}\{x,\log x\}$ by setting
\[
T_0x^n(\log x)^m=nx^{n-1}(\log x)^m\;\;\mbox{ and }\;\;
T_1x^n(\log x)^m=mx^{n-1}(\log x)^{m-1},
\]
respectively, and extending to all of ${\mathbb C}\{x,\log x\}$ by formal
linearity. Then the formula
\[
\frac d{dx}x^n(\log x)^m=nx^{n-1}(\log x)^m+mx^{n-1}(\log x)^{m-1}
\]
(extended to ${\mathbb C}\{x,\log x\}$) can be written as
\[
\frac d{dx}=T_0+T_1
\]
on ${\mathbb C}\{x,\log x\}$. So for $k\geq 1$, on ${\mathbb C}\{x,\log x\}$,
\begin{eqnarray*}
&&\bigg(\frac d{dx}\bigg)^k=\sum_{(i_1,\dots,i_k)\in\{0,1\}^k}
T_{i_1}\cdots T_{i_k}\\
&=&\sum_{j=0}^{k-1}\,\sum_{0\leq t_1<t_2<\cdots<t_{k-j}<k}
T_1^{k-t_{k-j}-1}T_0T_1^{t_{k-j}-t_{k-j-1}-1}T_0\cdots
T_0T_1^{t_2-t_1-1}T_0T_1^{t_1}+T_1^k,
\end{eqnarray*}
where $j$ gives the number of $T_1$'s in the product $T_{i_1}\cdots
T_{i_k}$; there are $k-j$ $T_0$'s, which are in the following
positions, reading {}from the right: $t_1+1$, $t_2+1, \dots,$
$t_{k-j}+1$. That is, for $k\geq 1$,
\begin{eqnarray*}
\lefteqn{\bigg(\frac d{dx}\bigg)^kx^n(\log x)^m=
\sum_{j=0}^km(m-1)\cdots(m-j+1)\cdot}\\
&&\cdot\bigg(\sum_{0\leq t_1<t_2<\cdots<t_{k-j}<k}
(n-t_1)(n-t_2)\cdots(n-t_{k-j})\bigg) x^{n-k}(\log x)^{m-j},
\end{eqnarray*}
where it is understood that if $j=k$, then the latter sum (in
parentheses) is $1$. In this formula, setting $n=0$, multiplying by
$y^k/k!$, and then summing over $k\in {\mathbb N}$, we get (noting that
$t_1=0$ contributes $0$)
\begin{eqnarray}\label{log:alg1}
\lefteqn{e^{y\frac d{dx}}(\log x)^m=\sum_{k\in {\mathbb N}}\bigg(\frac
yx\bigg)^k
\sum_{j=0}^k{m\choose j}(\log x)^{m-j}\frac{j!}{k!}\cdot}\nno\\
&&\cdot\bigg(\sum_{0<t_1<t_2<\cdots<t_{k-j}<k}
(-t_1)(-t_2)\cdots(-t_{k-j})\bigg).
\end{eqnarray}
So (\ref{log:ck1}) for $f(x)=(\log x)^m$ is equivalent to equating the
right-hand side of (\ref{log:alg1}) to
\begin{eqnarray}\label{log:alg2}
\lefteqn{\bigg(\log x+\sum_{i\geq 1} \frac{(-1)^{i-1}}i \bigg(\frac
  yx\bigg)^i\bigg)^m=}\nno\\
&&=\sum_{j\in {\mathbb N}}{m\choose j}(\log x)^{m-j} \bigg(\sum_{i\geq
1}\frac{(-1)^{i-1}}i \bigg(\frac yx\bigg)^i\bigg)^j\nno\\
&&=\sum_{j\in {\mathbb N}}{m\choose j}(\log x)^{m-j}\sum_{k\geq j}
\Bigg(\sum_{
\mbox{
\tiny
$\begin{array}{c}i_1+\cdots+i_j=k\\1\leq
i_1,\dots,i_j\leq k\end{array}$
}
} \frac{(-1)^{k-j}}{i_1i_2\cdots
i_j}\Bigg)\bigg(\frac yx\bigg)^k\nno\\
&&=\sum_{k\in {\mathbb N}}\bigg(\frac yx\bigg)^k\sum_{j=0}^k{m\choose
j}(\log x)^{m-j}\Bigg(\sum_{
\mbox{
\tiny
$\begin{array}{c}i_1+\cdots+i_j=k\\1
\leq i_1,\dots,i_j\leq k\end{array}$
}
} \frac{(-1)^{k-j}}{i_1i_2\cdots
i_j}\Bigg).
\end{eqnarray}
Comparing the right-hand sides of (\ref{log:alg1}) and
(\ref{log:alg2}) we see that it is equivalent to proving the
combinatorial identity
\begin{eqnarray}\label{log:comb}
\frac{j!}{k!}\sum_{0<t_1<t_2<\cdots<t_{k-j}<k} t_1t_2\cdots
t_{k-j}= \sum_{
\mbox{
\tiny
$\begin{array}{c}i_1+\cdots+i_j=k\\1\leq
i_1,\dots,i_j\leq k\end{array}$
}
} \frac{1}{i_1i_2\cdots i_j}
\end{eqnarray}
for all $k\in {\mathbb N}$ and $j=0,\dots,k$.
Note that there is no $m$ involved here.  But for $m$ a
nonnegative integer, (\ref{log:ck1}) for $f(x)=(\log x)^m$ follows
{}from
\[
e^{y\frac d{dx}}(\log x)^m=(e^{y\frac d{dx}}\log x)^m
\]
(recall Remark \ref{ddxchk}) and
\[
e^{y\frac d{dx}}\log x=\log x+\sum_{i\geq 1}\frac{(-1)^{i-1}}i \bigg(\frac
yx\bigg)^i
\]
(recall (\ref{logx+y})).  Thus the expressions in (\ref{log:alg1}) and
(\ref{log:alg2}) are equal for any such $m$.  Equating coefficients
and choosing $m\geq j$ gives us (\ref{log:comb}).  Therefore
(\ref{log:ck1}) also holds for any $m\in{\mathbb C}$.  \epf

\begin{rema}{\rm
Here is an amusing sidelight: When we were writing up the proof above,
one of us (L.Z.) happened to pick up the then-current issue of the American
Mathematical Monthly and happened to notice the following problem
{}from the Problems and Solutions section, proposed by D. Lubell
\cite{Lu}:
\begin{quote}
Let $N$ and $j$ be positive integers, and let $S=\{(w_1,\dots, w_j)\in
{\mathbb Z}_+^j\,|\,0<w_1+\cdots+w_j\leq N\}$ and $T=\{(w_1,\cdots,w_j)\in
{\mathbb Z}_+^j\,|\linebreak w_1,\dots,w_j \mbox{ are distinct and bounded by }N\}.$
Show that
$$\sum_S\frac 1{w_1\cdots w_j}=\sum_T\frac 1{w_1\cdots w_j}.$$
\end{quote}
But this follows immediately {}from (\ref{log:comb}) (which is in fact a
refinement), since the left-hand side of (\ref{log:comb}) is equal to
\[
j!\sum_{1\leq w_1<w_2<\cdots<w_{j-1}\leq k-1} \frac
1{w_1w_2\cdots w_{j-1}k}=\sum_{T_k} \frac 1{w_1w_2\cdots w_j}
\]
where
\[
T_k=\{(w_1,\dots,w_j)\in\{1,2,\dots,k\}^j\,|\,w_i\mbox{ distinct,
with maximum exactly }k\},
\]
the right-hand side is
\[
\sum_{S_k} \frac 1{w_1w_2\cdots w_j}
\]
where
\[
S_k=\{(w_1,\dots,w_j)\in\{1,2,\dots,k\}^j\,|\,
w_1+\cdots+w_j=k\},
\]
and one has $S=\coprod_{k=1}^N S_k$ and $T=\coprod_{k=1}^N T_k$. }
\end{rema}

When we define the notion of logarithmic intertwining operator below,
we will impose a condition requiring certain formal series to lie in
spaces of the type ${\cal W}[\log x]\{x\}$ (so that for each power of
$x$, possibly complex, we have a {\it polynomial} in $\log x$), partly
because such results as the following (which is expected) will indeed
hold in our formal setup when the powers of the formal variables are
restricted in this way (cf.\ Remark \ref{log:[[]]} below).

\begin{lemma}\label{log:de}
Let $a\in {\mathbb C}$ and $m\in {\mathbb Z}_+$. If $f(x)\in {\cal W}[\log
x]\{x\}$ (${\cal W}$ any vector space not involving $x$ or $\log x$)
satisfies the formal differential equation
\begin{equation}\label{de:(xdx-a)^m}
\bigg(x\frac{d}{dx}-a\bigg)^mf(x)=0,
\end{equation}
then $f(x)\in {\cal W}x^a\oplus{\cal W}x^a\log x \oplus\cdots\oplus{\cal
W}x^a(\log x)^{m-1}$, and furthermore, if $m$ is the smallest integer
so that (\ref{de:(xdx-a)^m}) is satisfied, then the coefficient of
$x^a(\log x)^{m-1}$ in $f(x)$ is nonzero.
\end{lemma}
\pf For any $f(x)=\sum_{n,k}w_{n,k}x^n(\log x)^k\in {\cal W}\{x,\log
x\}$,
\begin{eqnarray*}
x\frac d{dx}f(x)&=&\sum_{n,k}nw_{n,k}x^n(\log x)^k+\sum_{n,k}kw_{n,k}
x^n(\log x)^{k-1}\\
&=&\sum_{n,k}(nw_{n,k}+(k+1)w_{n,k+1})x^n(\log x)^k.
\end{eqnarray*}
Thus for any $a\in {\mathbb C}$,
\begin{equation}\label{de:act1}
\bigg(x\frac d{dx}-a\bigg)f(x)=\sum_{n,k}((n-a)w_{n,k}+(k+1)w_{n,k+1})
x^n(\log x)^k.
\end{equation}

Now suppose that $f(x)$ lies in ${\cal W}[\log x]\{x\}$. Let us prove
the assertion of the lemma by induction on $m$.

If $m=1$, by (\ref{de:act1}) we see that $(x\frac{d}{dx}-a)f(x)=0$
means that
\begin{equation}\label{de:m=1}
(n-a)w_{n,k}+(k+1)w_{n,k+1}=0\;\;\mbox{ for any }n\in{\mathbb C},\,
k\in{\mathbb Z}.
\end{equation}
Fix $n$. If $w_{n,k}\neq 0$ for some $k$, let $k_n$ be the smallest
nonnegative integer such that $w_{n,k}=0$ for any $k>k_n$ (such a
$k_n$ exists because $f(x)\in {\cal W}[\log x]\{x\}$). Then
\[
(n-a)w_{n,k_n}=-(k_n+1)w_{n,k_n+1}=0.
\]
But $w_{n,k_n}\neq 0$ by the choice of $k_n$, so we must have
$n=a$. Now (\ref{de:m=1}) becomes $(k+1)w_{a,k+1}=0$ for any $k\in
{\mathbb Z}$, so that $w_{a,k}=0$ unless $k=0$. Thus $f(x)=w_{a,0}x^a$.
If in addition $m=1$ is the smallest integer such that
(\ref{de:(xdx-a)^m}) holds, then $f(x)\neq 0$. So $w_{a,0}$, the
coefficient of $x^a$, is not zero.

Suppose the statement is true for $m$. Then for the case $m+1$, since
\begin{equation}\label{de:m+1}
0=\bigg(x\frac d{dx}-a\bigg)^{m+1}f(x)=\bigg(x\frac
d{dx}-a\bigg)^m\bigg(x\frac d{dx}-a\bigg)f(x)
\end{equation}
implies that
\begin{equation}\label{de:wbar}
\bigg(x\frac d{dx}-a\bigg)f(x)=\bar w_0x^a+\bar w_1x^a\log x+\cdots +
\bar w_{m-1}x^a(\log x)^{m-1}
\end{equation}
for some $\bar w_0, \bar w_1,\cdots, \bar w_{m-1}\in {\cal W}$, by
(\ref{de:act1}) we get
\begin{eqnarray*}
(n-a)w_{n,j}+(j+1)w_{n,j+1}=0&\;\;&\mbox{for any }n\neq a\mbox{ and
    any }j\in{\mathbb Z}\\
(j+1)w_{a,j+1}=\bar w_j&\;&\mbox{for any }j\in\{0,1,\dots,m-1\}\\
(j+1)w_{a,j+1}=0&\;&\mbox{for any }j\notin\{0,1,\dots,m-1\}
\end{eqnarray*}
By the same argument as above we get $w_{n,j}=0$ for any $n\neq a$ and
any $j$. So
\[
f(x)=w_{a,0}x^a+\bar w_0x^a\log x+\frac{\bar w_1}{2}x^a(\log x)^2+\cdots
+\frac{\bar w_{m-1}}mx^a(\log x)^m,
\]
as we want.  If in addition $m+1$ is the smallest integer so that
(\ref{de:m+1}) is satisfied, then by the induction assumption, $\bar
w_{m-1}$ in (\ref{de:wbar}) is not zero. So the coefficient in $f(x)$
of $x^a(\log x)^m$, $\bar w_{m-1}/m$, is not zero, as we want.  \epf

\begin{rema}\label{log:[[]]}{\rm
Note that there are solutions of the equation (\ref{de:(xdx-a)^m})
outside ${\cal W}[\log x]\{x\}$, for example,
\[
f(x)=wx^be^{(a-b)\log x}\in x^b{\cal W}[[\log x]]
\]
for any complex number $b\neq a$ and any
$0\neq w\in {\cal W}$.}
\end{rema}

Following \cite{Mi}, with a slight generalization (see Remark
\ref{log:compM}), we now introduce the notion of logarithmic
intertwining operator, together with the notion of
``grading-compatible logarithmic intertwining operator,'' adapted to
the strongly graded case.  We will later see that the axioms in these
definitions correspond exactly to those in the notion of certain
``intertwining maps'' (see Definition \ref{im:imdef} below).

\begin{defi}\label{log:def}{\rm
Let $(W_1,Y_1)$, $(W_2,Y_2)$ and $(W_3,Y_3)$ be generalized modules
for a M\"obius (or conformal) vertex algebra $V$. A {\em logarithmic
intertwining operator of type ${W_3\choose W_1\,W_2}$} is a linear map
\begin{equation}\label{log:map0}
{\cal Y}(\cdot, x)\cdot: W_1\otimes W_2\to W_3[\log x]\{x\},
\end{equation}
or equivalently,
\begin{equation}\label{log:map}
w_{(1)}\otimes w_{(2)}\mapsto{\cal Y}(w_{(1)},x)w_{(2)}=\sum_{n\in {\mathbb
C}}\sum_{k\in {\mathbb N}}{w_{(1)}}_{n;\,k}^{\cal Y}w_{(2)}x^{-n-1}(\log
x)^k\in W_3[\log x]\{x\}
\end{equation}
for all $w_{(1)}\in W_1$ and $w_{(2)}\in W_2$, such that the following
conditions are satisfied: the {\em lower truncation condition}: for
any $w_{(1)}\in W_1$, $w_{(2)}\in W_2$ and $n\in {\mathbb C}$,
\begin{equation}\label{log:ltc}
{w_{(1)}}_{n+m;\,k}^{\cal Y}w_{(2)}=0\;\;\mbox{ for }\;m\in {\mathbb N}
\;\mbox{ sufficiently large,\, independently of}\;k;
\end{equation}
the {\em Jacobi identity}:
\begin{eqnarray}\label{log:jacobi}
\lefteqn{\dps x^{-1}_0\delta \bigg( {x_1-x_2\over x_0}\bigg)
Y_3(v,x_1){\cal Y}(w_{(1)},x_2)w_{(2)}}\nno\\
&&\hspace{2em}- x^{-1}_0\delta \bigg( {x_2-x_1\over -x_0}\bigg)
{\cal Y}(w_{(1)},x_2)Y_2(v,x_1)w_{(2)}\nno \\
&&{\dps = x^{-1}_2\delta \bigg( {x_1-x_0\over x_2}\bigg) {\cal
Y}(Y_1(v,x_0)w_{(1)},x_2) w_{(2)}}
\end{eqnarray}
for $v\in V$, $w_{(1)}\in W_1$ and $w_{(2)}\in W_2$ (note that the
first term on the left-hand side is meaningful because of
(\ref{log:ltc})); the {\em $L(-1)$-derivative property:} for any
$w_{(1)}\in W_1$,
\begin{equation}\label{log:L(-1)dev}
{\cal Y}(L(-1)w_{(1)},x)=\frac d{dx}{\cal Y}(w_{(1)},x);
\end{equation}
and the {\em ${\mathfrak s}{\mathfrak l}(2)$-bracket relations:} for any
$w_{(1)}\in W_1$,
\begin{equation}\label{log:L(j)b}
{}[L(j), {\cal Y}(w_{(1)},x)]=\sum_{i=0}^{j+1}{j+1\choose i}x^i{\cal
Y}(L(j-i)w_{(1)},x)
\end{equation}
for $j=-1, 0$ and $1$ (note that if $V$ is in fact a conformal vertex
algebra, this follows automatically {}from the Jacobi identity
(\ref{log:jacobi}) by setting $v=\omega$ and then taking
$\res_{x_0}\res_{x_1}x_1^{j+1}$).  }
\end{defi}

\begin{rema}{\rm
We will sometimes write the Jacobi identity (\ref{log:jacobi}) as
\begin{eqnarray}
\lefteqn{\dps x^{-1}_0\delta \bigg( {x_1-x_2\over x_0}\bigg)
Y(v,x_1){\cal Y}(w_{(1)},x_2)w_{(2)}}\nno\\
&&\hspace{2em}- x^{-1}_0\delta \bigg( {x_2-x_1\over -x_0}\bigg)
{\cal Y}(w_{(1)},x_2)Y(v,x_1)w_{(2)}\nno \\
&&{\dps = x^{-1}_2\delta \bigg( {x_1-x_0\over x_2}\bigg) {\cal
Y}(Y(v,x_0)w_{(1)},x_2) w_{(2)}}
\end{eqnarray}
(dropping the subscripts on the module actions) for brevity.}
\end{rema}

\begin{rema}\label{ordinaryandlogintwops}{\rm
The ordinary intertwining operators (as in, for example,
\cite{tensor1}) among triples of modules for a vertex operator algebra
are exactly the logarithmic intertwining operators that do not involve
the formal variable $\log x$, except for our present relaxation of the
lower truncation condition.  The lower truncation condition that we
use here can be equivalently stated as: For any $w_{(1)}\in W_1$,
$w_{(2)}\in W_2$ and $n\in{\mathbb C}$, there is no nonzero term
involving $x^{n-m}$ appearing in ${\cal Y}(w_{(1)},x)w_{(2)}$ when
$m\in{\mathbb N}$ is large enough.  In \cite{tensor1}, the lower
truncation condition in the definition of the notion of intertwining
operator states: For any $w_{(1)}\in W_1$ and $w_{(2)}\in W_2$,
\begin{equation}\label{0forlargen}
(w_{(1)})_nw_{(2)}=0\;\mbox{for}\;n\;\mbox{whose real part is
sufficiently large}.
\end{equation}
This is slightly stronger than the lower truncation condition that we
use here, even if no $\log x$ is involved, when the powers of $x$ in
${\cal Y}(w_{(1)},x)w_{(2)}$ belong to infinitely many different
congruence classes modulo ${\mathbb Z}$.  But in the case that our
(generalized) modules satisfy the lower boundedness condition
(\ref{set:dmltc-1}), this stronger condition (\ref{0forlargen}) will
hold.}
\end{rema}

\begin{rema}\label{g-mod-as-l-int}{\rm
Given a generalized module $(W, Y_{W})$ for a M\"{o}bius (or conformal) 
vertex algebra $V$, the vertex operator map $Y_{W}$ itself is clearly 
a logarithmic intertwining operator of type ${W\choose VW}$; in fact, it 
does not involve $\log x$ and its powers of $x$ are all integers. In 
particular, taking $(W, Y_{W})$ to be $(V, Y)$ itself, we have that the 
vertex operator map $Y$ is a logarithmic intertwining operator of type
${V\choose VV}$ not involving $\log x$ and having only integral powers 
of $x$.}
\end{rema}

The logarithmic intertwining operators of a fixed type ${W_3\choose
W_1\,W_2}$ form a vector space.

\begin{defi}\label{gradingcompatintwop}{\rm
In the setting of Definition \ref{log:def}, suppose in addition that
$V$ and $W_1$, $W_2$ and $W_3$ are strongly graded (recall Definitions
\ref{def:dgv} and \ref{def:dgw}).  A logarithmic intertwining operator
${\cal Y}$ as in Definition \ref{log:def} is a {\em grading-compatible
logarithmic intertwining operator} if for $\beta, \gamma \in \tilde A$
(recall Definition \ref{def:dgw}) and $w_{(1)} \in W_{1}^{(\beta)}$,
$w_{(2)} \in W_{2}^{(\gamma)}$, $n \in \C$ and $k \in \mathbb N$, we
have
\begin{equation}\label{gradingcompatcondn}
{w_{(1)}}_{n;\,k}^{\cal Y}w_{(2)} \in W_{3}^{(\beta + \gamma)}.
\end{equation}}
\end{defi}

\begin{rema}{\rm
The term ``grading-compatible'' in Definition
\ref{gradingcompatintwop} refers to the $\tilde A$-gradings; {\it any}
logarithmic intertwining operator is compatible with the $\C$-gradings
of $W_1$, $W_2$ and $W_3$, in view of Proposition \ref{log:logwt}(b)
below.}
\end{rema}

\begin{rema}\label{str-graded-g-mod-as-l-int}
{\rm Given a strongly graded generalized module $(W, Y_{W})$ for a
strongly graded M\"{o}bius (or conformal) vertex algebra $V$, the
vertex operator map $Y_{W}$ is a grading-compatible logarithmic
intertwining operator of type ${W\choose VW}$ not involving $\log x$
and having only integral powers of $x$.  Taking $(W, Y_{W})$ in
particular to be $(V, Y)$ itself, we have that the vertex operator map
$Y$ is a grading-compatible logarithmic intertwining operator of type
${V\choose VV}$ not involving $\log x$ and having only integral powers
of $x$.}
\end{rema}

In the strongly graded context (the main context for our tensor
product theory), we will use the following notation and terminology,
traditionally used in the setting of ordinary intertwining operators,
as in \cite{FHL}:

\begin{defi}\label{fusionrule}{\rm
In the setting of Definition \ref{gradingcompatintwop}, the
grading-compatible logarithmic intertwining operators of a fixed type
${W_3\choose W_1\, W_2}$ form a vector space, which we denote by
${\cal V}^{W_3}_{W_1\,W_2}$.  We call the dimension of ${\cal
V}^{W_3}_{W_1\,W_2}$ the {\it fusion rule} for $W_1$, $W_2$ and $W_3$
and denote it by $N^{W_3}_{W_1\,W_2}$.}
\end{defi}

\begin{rema}{\rm
In the strongly graded context, suppose that $W_1$, $W_2$ and $W_3$ in
Definition \ref{log:def} are expressed as finite direct sums of
submodules.  Then the space ${\cal V}^{W_3}_{W_1\,W_2}$ can be
naturally expressed as the corresponding (finite) direct sum of the
spaces of (grading-compatible) logarithmic intertwining operators
among the direct summands, and the fusion rule $N^{W_3}_{W_1\,W_2}$ is
thus the sum of the fusion rules for the direct summands.
}
\end{rema}

\begin{rema}{\rm
As we shall point out in Remark \ref{log:ordi} below, it turns out
that the notion of fusion rule in Definition \ref{fusionrule} agrees
with the traditional notion, in the case of a vertex operator algebra
and ordinary modules.  The justification of this assertion uses Parts
(b) and (c), or alternatively, Part (a), of the next proposition.
Part (a), whose proof uses Lemma \ref{log:de}, shows how logarithmic
intertwining operators yield expansions involving only finitely many
powers of $\log x$.  Part (b) is a generalization of formula
(\ref{set:wtvn}).
}
\end{rema}

\begin{propo}\label{log:logwt}
Let $W_1$, $W_2$, $W_3$ be generalized modules for a M\"obius (or
conformal) vertex algebra $V$, and let ${\cal Y}(\cdot, x)\cdot$ be a
logarithmic intertwining operator of type ${W_3\choose W_1\,W_2}$.
Let $w_{(1)}$ and $w_{(2)}$ be homogeneous elements of $W_1$ and $W_2$
of generalized weights $n_1$ and $n_2 \in {\mathbb C}$, respectively, and
let $k_1$ and $k_2$ be positive integers such that
\[
(L(0)-n_1)^{k_1}w_{(1)}=0 \;\;\mbox{and} \;\; (L(0)-n_2)^{k_2}w_{(2)}=0.
\]
Then we have:

(a) (\cite{Mi}) For any $w'_{(3)}\in W_3^*$, $n_3\in {\mathbb C}$ and
$k_3\in {\mathbb Z}_+$ such that $(L'(0)-n_3)^{k_3}w'_{(3)}=0$,
\begin{eqnarray}\label{log:k}
\lefteqn{\langle w'_{(3)}, {\cal Y}(w_{(1)}, x)w_{(2)}\rangle}\nno\\
&&\in {\mathbb C}x^{n_3-n_1-n_2}\oplus{\mathbb C}x^{n_3-n_1-n_2}\log
x\oplus\cdots\oplus {\mathbb C}x^{n_3-n_1-n_2}(\log
x)^{k_1+k_2+k_3-3}.\nno\\
\end{eqnarray}

(b) For any $n\in {\mathbb C}$ and $k\in {\mathbb N}$, ${w_{(1)}}^{\cal
Y}_{n;\,k}w_{(2)}\in W_3$ is homogeneous of generalized weight
$n_1+n_2-n-1$.

(c) Fix $n\in {\mathbb C}$ and $k\in {\mathbb N}$. For each $i,j\in {\mathbb
N}$, let $m_{ij}$ be a nonnegative integer such that
\[
(L(0)-n_1-n_2+n+1)^{m_{ij}}(((L(0)-n_1)^iw_{(1)})^{\cal
Y}_{n;\,k}(L(0)-n_2)^jw_{(2)})=0.
\]
Then for all $t\geq \max\{m_{ij}\,|\,0\leq i<k_1,\; 0\leq
j<k_2\}+k_1+k_2-2$,
\[
{w_{(1)}}^{\cal Y}_{n;\,k+t}w_{(2)}=0.
\]
\end{propo}

We will need the following lemma in the proof:

\begin{lemma}\label{log:lemma}
Let $W_1$, $W_2$, $W_3$ be generalized modules for a M\"obius (or
conformal) vertex algebra $V$. Let
\begin{eqnarray}
{\cal Y}(\cdot,x)\cdot: W_1\otimes
W_2 & \to & W_3\{x,\log x\} \nonumber\\
w_{(1)}\otimes w_{(2)} & \mapsto & {\cal Y}(w_{(1)},x)w_{(2)}=\sum_{n,k\in
{\mathbb C}}{w_{(1)}}_{n;\,k}^{\cal Y}w_{(2)}x^{-n-1}(\log x)^k
\nonumber\\
& & \label{intertwopinlemma}
\end{eqnarray}
be a linear map that satisfies the L(-1)-derivative property
(\ref{log:L(-1)dev}) and the $L(0)$-bracket relation, that is,
(\ref{log:L(j)b}) with $j=0$. Then for any $a,b,c\in{\mathbb C}$,
$t\in{\mathbb N}$, $w_{(1)}\in W_1$ and $w_{(2)}\in W_2$,
\begin{eqnarray}\label{log:ty}
\lefteqn{(L(0)-c)^t{\cal Y}(w_{(1)},x)w_{(2)}=\sum_{i,j,l\in {\mathbb
N},\;i+j+l=t}\frac{t!}{i!j!l!}\cdot}\nno\\
&&\cdot\left(x\frac d{dx}-c+a+b\right)^l{\cal Y}((L(0)-a)^iw_{(1)},x)
(L(0)-b)^jw_{(2)}.\nno\\
\end{eqnarray}
Also, for any $a,b,n,k\in{\mathbb C}$, $t\in{\mathbb N}$,
$w_{(1)}\in W_1$ and $w_{(2)}\in W_2$, we have
\begin{eqnarray}\label{log:t00}
\lefteqn{(L(0)-a-b+n+1)^t({w_{(1)}}_{n;\,k}^{\cal Y}w_{(2)})}\nno\\
&&=t!\sum_{i,j,l\geq 0,\;i+j+l=t}{k+l\choose l}
\bigg(\frac{(L(0)-a)^i}{i!}w_{(1)}\bigg)^{\cal
Y}_{n;\,k+l}\bigg(\frac{(L(0)-b)^j}{j!}w_{(2)}\bigg);\nno\\
\end{eqnarray}
in generating function form, this gives
\begin{eqnarray}\label{log:e^L(0)}
\lefteqn{e^{y(L(0)-a-b+n+1)}({w_{(1)}}_{n;\,k}^{\cal
Y}w_{(2)})}\nno\\ 
&&=\sum_{l\in{\mathbb N}}{k+l\choose
l}(e^{y(L(0)-a)}w_{(1)})_{n;\,k+l}^{\cal Y}(e^{y(L(0)-b)}w_{(2)})y^l.
\end{eqnarray}
\end{lemma}
\pf {}From (\ref{log:L(-1)dev}) and (\ref{log:L(j)b}) with $j=0$ we have
\begin{eqnarray}\label{log:pf1}
\lefteqn{L(0){\cal Y}(w_{(1)},x)w_{(2)}={\cal Y}(w_{(1)},x)
L(0)w_{(2)}}\nno\\
&&+x\frac d{dx}{\cal Y}(w_{(1)},x)w_{(2)}+
{\cal Y}(L(0)w_{(1)},x)w_{(2)}.
\end{eqnarray}
Hence
\begin{eqnarray*}
\lefteqn{(L(0)-c){\cal Y}(w_{(1)},x)w_{(2)}={\cal Y}(w_{(1)},x)
(L(0)-b)w_{(2)}}\\
&&+\left(x\frac d{dx}-c+a+b\right){\cal Y}(w_{(1)},x)w_{(2)}+ {\cal
Y}((L(0)-a)w_{(1)},x)w_{(2)}
\end{eqnarray*}
for any complex numbers $a$, $b$ and $c$. In view of the fact that the
actions of $L(0)$ and $d/dx$ commute with each other, this implies
(\ref{log:ty}) essentially because of the expansion formula for powers
of a sum of commuting operators, that is, for any commuting operators
$T_1, \dots, T_s$ and $t\in {\mathbb N}$,
\begin{equation}\label{log:expand}
(T_1+\cdots+T_s)^t=\sum_{i_1,\dots,i_s\in{\mathbb N},\; i_1+\cdots+i_s=t}
\frac{t!}{i_1!\cdots i_s!}T_1^{i_1}\cdots T_s^{i_s}.
\end{equation}

On the other hand, by taking coefficient of $x^{-n-1}(\log x)^k$ in
(\ref{log:pf1}) we get
\begin{eqnarray*}
\lefteqn{L(0){w_{(1)}}_{n;\,k}^{\cal Y}w_{(2)}={w_{(1)}}_{n;\,k}^{\cal
Y}L(0)w_{(2)}+(-n-1){w_{(1)}}_{n;\,k}^{\cal Y}w_{(2)}}\nno\\
&&+(k+1){w_{(1)}}_{n;\,k+1}^{\cal
Y}w_{(2)}+(L(0)w_{(1)})_{n;\,k}^{\cal Y}w_{(2)}.
\end{eqnarray*}
So for any $a,b,n,k\in {\mathbb C}$,
\begin{eqnarray}\label{log:t=1}
\lefteqn{(L(0)-a-b+n+1)({w_{(1)}}_{n;\,k}^{\cal Y}w_{(2)})=
((L(0)-a)w_{(1)})_{n;\,k}^{\cal Y}w_{(2)}}\nno\\
&&\hspace{2cm}+{w_{(1)}}_{n;\,k}^{\cal
Y}(L(0)-b)w_{(2)}+(k+1){w_{(1)}}_{n;\,k+1}^{\cal Y}w_{(2)}.
\end{eqnarray}
For $p,q\in {\mathbb N}$ and $n,k\in {\mathbb C}$, let us write
\begin{equation}\label{log:Tprq}
T_{p,k,q}=((L(0)-a)^pw_{(1)})^{\cal Y}_{n;\,k}((L(0)-b)^qw_{(2)}).
\end{equation}
Then {}from (\ref{log:t=1}) we see that for any $p,q\in {\mathbb N}$ and
$a,b,n,k\in {\mathbb C}$,
\begin{equation}\label{log:L(0)^1}
(L(0)-a-b+n+1)T_{p,k,q}=T_{p+1,k,q}+(k+1)T_{p,k+1,q}+T_{p,k,q+1}.
\end{equation}
Hence by (\ref{log:expand}) we have
\[
(L(0)-a-b+n+1)^tT_{p,k,q}=t!\sum_{i,j,l\geq 0,\;i+j+l=t}\frac
{(k+1)(k+2)\cdots(k+l)}{i!j!l!}T_{p+i,k+l,q+j}
\]
for any $a,b,n,k\in {\mathbb C}$ and $p,q\in {\mathbb N}$. In particular, by
setting $p=q=0$ we get (\ref{log:t00}), and (\ref{log:e^L(0)}) follows
easily {}from (\ref{log:t00}) by multiplying by $y^t/t!$ and then
summing over $t\in{\mathbb N}$.  \epfv

{\it Proof of Proposition \ref{log:logwt}}\hspace{2ex} (a): Under the
assumptions of the proposition, let us show that
\begin{equation}\label{log:ode}
\left\langle w'_{(3)},\left(x\frac d{dx}-n_3+n_1+n_2\right)^
{k_3+k_1+k_2-2}{\cal Y}(w_{(1)},x)w_{(2)}\right\rangle=0
\end{equation}
by induction on $k_1+k_2$.

For $k_1=k_2=1$, {}from (\ref{log:ty}) with $a=n_1$, $b=n_2$, $c=n_3$
and $t=k_3$ we have
\begin{eqnarray*}
0&=&\langle (L'(0)-n_3)^{k_3}w'_{(3)},{\cal
Y}(w_{(1)},x)w_{(2)}\rangle\\
&=&\langle w'_{(3)},(L(0)-n_3)^{k_3}{\cal
Y}(w_{(1)},x)w_{(2)}\rangle\\
&=&\left\langle w'_{(3)},\left(x\frac d{dx}-n_3+n_1+n_2\right)^{k_3}{\cal
Y}(w_{(1)},x)w_{(2)}\right\rangle,
\end{eqnarray*}
which is (\ref{log:ode}) in the case $k_1=k_2=1$.

Suppose that (\ref{log:ode}) is true for all the cases with smaller
$k_1+k_2$. Then {}from (\ref{log:ty}) with $a=n_1$, $b=n_2$, $c=n_3$ and
$t=k_3+k_1+k_2-2$ we have
\begin{eqnarray*}
0&=&\langle (L'(0)-n_3)^{k_3+k_1+k_2-2}w'_{(3)},{\cal
Y}(w_{(1)},x)w_{(2)}\rangle\\
&=&\langle w'_{(3)},(L(0)-n_3)^{k_3+k_1+k_2-2}{\cal
Y}(w_{(1)},x)w_{(2)}\rangle\\
&=&\Biggl\langle w'_{(3)},\sum_{i,j,k\in {\mathbb
N},\; i+j+k=k_3+k_1+k_2-2}\frac{(k_3+k_1+k_2-2)!}{i!j!k!}\cdot\\
&&\quad\cdot\left(x\frac d{dx}-n_3+n_1+n_2\right)^k{\cal Y}
((L(0)-n_1)^iw_{(1)},x) (L(0)-n_2)^jw_{(2)}\Biggr\rangle\\
&=&\left\langle w'_{(3)},\left(x\frac d{dx}-n_3+n_1+n_2\right)^
{k_3+k_1+k_2-2}{\cal Y}(w_{(1)},x)w_{(2)}\right\rangle,
\end{eqnarray*}
where the last equality uses the induction assumption for the pair of
elements $(L(0)-n_1)^iw_{(1)}$ and $(L(0)-n_2)^jw_{(2)}$ for all
$(i,j)\neq (0,0)$. So (\ref{log:ode}) is established, that is, we have
the formal differential equation
\[
\left(x\frac d{dx}-n_3+n_1+n_2\right)^{k_3+k_1+k_2-2}\langle w'_{(3)},
{\cal Y}(w_{(1)},x)w_{(2)}\rangle=0.
\]
This implies (a) by Lemma \ref{log:de}.

(b): This follows {}from (\ref{log:t00}) with $a=n_1$, $b=n_2$ and the
fact that for any $\bar w_{(1)}\in W_1$, $\bar w_{(2)}\in W_2$ and
$\bar n\in {\mathbb C}$, there exists $K\in{\mathbb N}$ so that $(\bar
w_{(1)})^{\cal Y}_{\bar n;\bar k}\bar w_{(2)}=0$ for all $\bar k>K$,
due to (\ref{log:map0}).

(c): Let us prove (c) by induction on $k_1+k_2$ again. For
$k_1=k_2=1$, (\ref{log:t00}) with $a=n_1$ and $b=n_2$ gives
\[
(L(0)-n_1-n_2+n+1)^t({w_{(1)}}_{n;\,k}^{\cal
Y}w_{(2)})=((k+t)!/k!){w_{(1)}}_{n;\,k+t}^{\cal Y}w_{(2)},
\]
that is,
\[
{w_{(1)}}_{n;\,k+t}^{\cal Y}w_{(2)}= (k!/(k+t)!)
(L(0)-n_1-n_2+n+1)^t({w_{(1)}}_{n;\,k}^{\cal Y}w_{(2)}).
\]
So for $t\geq m_{00}$, ${w_{(1)}}_{n;\,k+t}^{\cal Y}w_{(2)}=0$,
proving the statement in case $k_1=k_2=1$.

Suppose that the statement (c) is true for all smaller $k_1+k_2$.
Then for $(i,j)\neq (0,0)$,
\[
\bigg(\frac{(L(0)-n_1)^i}{i!}w_{(1)}\bigg)^{\cal
Y}_{n;\,k+l}\bigg(\frac{(L(0)-n_2)^j}{j!}w_{(2)}\bigg)=0
\]
when $l\geq \max\{m_{i'j'}\,|\,i\leq i'<k_1,\; j\leq
j'<k_2\}+(k_1-i)+(k_2-j)-2$, and in particular, when $l\geq
\max\{m_{i'j'}\,|\,0\leq i'<k_1,\; 0\leq j'<k_2\}+k_1+k_2-i-j-2$. But
then for all $t\geq\max\{m_{ij}\,|\,0\leq i<k_1,\; 0\leq
j<k_2\}+k_1+k_2-2$, (\ref{log:t00}) gives
$0=((k+t)!/k!){w_{(1)}}_{n;\,k+t}^{\cal Y}w_{(2)}$, proving what we
need.  \epfv

The following corollary is immediate {}from Proposition \ref{log:logwt}(b):

\begin{corol}\label{powerscongruentmodZ}
Let $V$ be a M\"obius (or conformal) vertex algebra and let $W_1$,
$W_2$ and $W_3$ be generalized $V$-modules whose weights are all
congruent modulo ${\mathbb Z}$ to complex numbers $h_1$, $h_2$ and $h_3$,
respectively.  (For example, $W_1$, $W_2$ and $W_3$ might be
indecomposable; recall Remark \ref{congruent}.)  Let ${\cal Y}(\cdot,
x)\cdot$ be a logarithmic intertwining operator of type ${W_3\choose
W_1\,W_2}$.  Then all powers of $x$ in ${\cal Y}(\cdot, x)\cdot$ are
congruent modulo ${\mathbb Z}$ to $h_3-h_1-h_2$. \epf
\end{corol}

\begin{rema}\label{log:ordi}{\rm
Let $W_1$, $W_2$ and $W_3$ be (ordinary) modules for a M\"obius (or
conformal) vertex algebra $V$. Then any logarithmic intertwining
operator of type ${W_3\choose W_1\,W_2}$ is just an ordinary
intertwining operator of this type, i.e., it does not involve $\log
x$. This clearly follows {}from Proposition \ref{log:logwt}(b) and (c),
where $k_1$ and $k_2$ are chosen to be $1$, $k$ is chosen to be $0$,
and $m_{00}$ is chosen to be $1$.  It also follows, alternatively,
{}from Proposition \ref{log:logwt}(a).  As a result, for $V$ a vertex
operator algebra (viewed as a conformal vertex algebra strongly graded
with respect to the trivial group; recall Remark \ref{rm1}) and $W_1$,
$W_2$ and $W_3$ $V$-modules in the sense of Remark
\ref{moduleswiththetrivialgroup}, the notion of fusion rule defined in
this work (recall Definition \ref{fusionrule}) coincides with the
notion of fusion rule defined in, for example, \cite{tensor1} (except
for the minor issue of the truncation condition for an intertwining
operator, discussed in Remark \ref{ordinaryandlogintwops}).}
\end{rema}

\begin{rema}\label{log:compM}{\rm
Our definition of logarithmic intertwining operator is identical to
that in \cite{Mi} (in case $V$ is a vertex operator algebra) except
that in \cite{Mi}, a logarithmic intertwining operator ${\cal
Y}$ of type ${W_3\choose W_1\,W_2}$ is required to be a
linear map $W_1\to \hom(W_2,W_3)\{x\}[\log x]$, instead of as in
(\ref{log:map0}), and the lower truncation condition (\ref{log:ltc})
is replaced by: For any $w_{(1)}\in W_1$, $w_{(2)}\in W_2$ and $k\in
{\mathbb N}$,
\begin{equation}\label{repartbounded}
{w_{(1)}}^{\cal Y}_{n;k}w_{(2)}=0\;\;\mbox{for}\;n\;\mbox{whose real
part is sufficiently large}.
\end{equation}
Given generalized $V$-modules $W_1$, $W_2$ and $W_3$, suppose that for
each $i=1,2,3$, there exists some $K_i\in {\mathbb Z}_+$ such that
$(L(0)-L(0)_s)^{K_i}W_i=0$ (this is satisfied by many interesting
examples and is assumed in \cite{Mi} for all generalized modules under
consideration). Then for any logarithmic intertwining operator ${\cal
Y}$, any homogeneous elements $w_{(1)}\in {W_1}_{[n_1]}$, $w_{(2)}\in
{W_2}_{[n_2]}$, $n_1, n_2\in{\mathbb C}$, any $n\in{\mathbb C}$ and
any $k\in{\mathbb N}$, all the $m_{ij}$'s in Proposition
\ref{log:logwt}(c) can be chosen to be no greater than $K_3$, while
$k_1$ and $k_2$ can be chosen to be no greater than $K_1$ and $K_2$,
respectively.  Proposition \ref{log:logwt}(c) thus
implies that the largest power of $\log x$ that is involved in ${\cal
Y}$ is no greater than $K_1+K_2+K_3-3$.  (This follows alternatively
{}from Proposition \ref{log:logwt}(a) and (\ref{L(0)N}).)
In particular, ${\cal Y}$ maps
$W_1$ to $\hom(W_2,W_3)\{x\}[\log x]$, and in fact, we even have that
$K_1+K_2+K_3-3$ is a global bound on the powers of $\log x$,
independently of $w_{(1)}\in W_1$, so that
\begin{equation}
{\cal Y}(\cdot,x)\cdot \in \hom(W_1 \otimes W_2,W_3)\{x\}[\log x],
\end{equation}
and this global bound on the powers of $\log x$ is even independent of
${\cal Y}$.}
\end{rema}

\begin{rema}\label{=0}{\rm
In the setting of Definition \ref{gradingcompatintwop}, if $W_3$ is
lower bounded, then any grading-compatible logarithmic intertwining
operator ${\cal Y}$ of type ${W_3\choose W_1\,W_2}$ satisfies
(\ref{repartbounded}), in view of Proposition \ref{log:logwt}(b).}
\end{rema}

\begin{rema}\label{Y(k)}{\rm
Given a logarithmic intertwining operator ${\cal Y}$ as in
(\ref{log:map}), set
\[
{\cal Y}^{(k)}(w_{(1)},x)w_{(2)}=\sum_{n\in {\mathbb
C}}{w_{(1)}}_{n;\,k}^{\cal Y}w_{(2)}x^{-n-1}
\]
for $k\in {\mathbb N}$, $w_{(1)}\in W_1$ and $w_{(2)}\in W_2$, so that
\[
{\cal Y}(w_{(1)},x)w_{(2)}=\sum_{k\in {\mathbb N}}{\cal
  Y}^{(k)}(w_{(1)},x)w_{(2)}(\log x)^k.
\]
By taking the coefficients of the powers of $\log x_2$ and $\log x$ in
(\ref{log:jacobi}) and (\ref{log:L(j)b}), respectively, we see that
each ${\cal Y}^{(k)}$ satisfies the Jacobi identity and the ${\mathfrak
s}{\mathfrak l}(2)$-bracket relations. On the other hand, taking the
coefficients of the powers of $\log x$ in (\ref{log:L(-1)dev}) gives
\begin{equation}\label{log:L(-1)comp}
{\cal Y}^{(k)}(L(-1)w_{(1)},x)=\frac d{dx}{\cal Y}^{(k)}(w_{(1)},x)
+\frac{k+1}x{\cal Y}^{(k+1)}(w_{(1)},x)
\end{equation}
for any $k\in {\mathbb N}$ and $w_{(1)}\in W_1$. So ${\cal Y}^{(k)}$
does not in general satisfy the $L(-1)$-derivative property.  (If
${\cal Y}^{(k+1)}=0$, then ${\cal Y}^{(k)}$ of course does satisfy the
$L(-1)$-derivative property and so is an (ordinary) intertwining
operator; this certainly happens for $k=0$ in the context of Remark
\ref{log:ordi} and for $k=K_1+K_2+K_3-3$ in the context of Remark
\ref{log:compM}.)  However, in the following we will see that suitable
formal linear combinations of certain modifications of ${\cal
Y}^{(k)}$ (depending on $t\in{\mathbb N}$; see below) form a sequence
of logarithmic intertwining operators.  }
\end{rema}

\begin{rema}\label{log:mu}{\rm
Given a logarithmic intertwining operator ${\cal Y}$, let us write
\begin{eqnarray*}
{\cal Y}(w_{(1)},x)w_{(2)}&=&\sum_{n\in {\mathbb C}}\sum_{k\in {\mathbb N}}
{w_{(1)}}^{\cal Y}_{n;\,k}w_{(2)}x^{-n-1}(\log x)^k\\
&=&\sum_{\mu\in {\mathbb C}/{\mathbb Z}}\sum_{\bar n=\mu}\sum_{k\in {\mathbb
N}}{w_{(1)}}^{\cal Y}_{n;\,k}w_{(2)}x^{-n-1}(\log x)^k
\end{eqnarray*}
for any $w_{(1)}\in W_1$ and $w_{(2)}\in W_2$, where $\bar n$ denotes
the equivalence class of $n$ in ${\mathbb C}/{\mathbb Z}$. By extracting
summands corresponding to the same congruence class modulo ${\mathbb Z}$
of the powers of $x$ in (\ref{log:jacobi}), (\ref{log:L(-1)dev}) and
(\ref{log:L(j)b}) we see that for each $\mu\in {\mathbb C}/{\mathbb Z}$,
\begin{equation}\label{log:c+n}
w_{(1)}\otimes w_{(2)}\mapsto {{\cal Y}^\mu}(w_{(1)},x)w_{(2)}= 
\sum_{\bar n=\mu}\sum_{k\in {\mathbb N}}
{w_{(1)}}^{\cal Y}_{n;\,k}w_{(2)}x^{-n-1}(\log x)^k
\end{equation}
still defines a logarithmic intertwining operator.  In the strongly
graded case, if ${\cal Y}$ is grading-compatible, then so is the
operator ${\cal Y}^\mu$ in (\ref{log:c+n}).  Conversely, suppose that
we are given a family of logarithmic intertwining operators $\{{\cal
Y}^\mu | \mu\in{\mathbb C}/{\mathbb Z}\}$ parametrized by
$\mu\in{\mathbb C}/{\mathbb Z}$ such that the powers of $x$ in ${\cal
Y}^\mu$ are restricted as in (\ref{log:c+n}). Then the formal sum
$\sum_{\mu\in{\mathbb C}/{\mathbb Z}}{\cal Y}^\mu$ is well defined and
is a logarithmic intertwining operator.  In the strongly graded case,
if each ${\cal Y}^\mu$ is grading-compatible, then so is this sum.}
\end{rema}

In the setting of Definition \ref{log:def}, for any integer $p$, set
\begin{equation}\label{substitutionofe2piipx}
{\cal Y} (w_{(1)},e^{2\pi ip}x)w_{(2)} = {\cal
Y}(w_{(1)},y)w_{(2)}\lbar_{y^n=e^{2\pi ipn}x^n,\; (\log y)^k=({2\pi
ip} +\log x)^k,\;n\in{\mathbb C},\;k\in{\mathbb N}}.
\end{equation}
This is in fact a well-defined element of $W_3[\log x]\{x\}$.  Note
that this element certainly depends on $p$, not just on $e^{2\pi ip}$
($=1$).  This substitution, which can be thought of as ``$x \mapsto
e^{2\pi ip}x$,'' will be considered in a more general form in
(\ref{log:subs}) below.

\begin{rema}\label{formalinvariance}{\rm
It is clear that in Definition \ref{log:def}, for any integer $p$, all
the axioms are formally invariant under the substitution
\[
x \mapsto e^{2\pi ip}x
\]
given by (\ref{substitutionofe2piipx}).  That is, if we
apply this substitution to each axiom, the axiom keeps the same form,
with the operator ${\cal Y}(\cdot,x)$ replaced by
\[
{\cal Y}(\cdot,e^{2\pi ip}x).
\]
For example, for the Jacobi identity
(\ref{log:jacobi}), we perform the substitution $x_{2} \mapsto e^{2\pi
ip}x_{2}$; the formal delta-functions remain unchanged because they
involve only integral powers of $x_2$ and no logarithms.  It follows
that ${\cal Y}(\cdot,e^{2\pi ip}x)$ is again a logarithmic
intertwining operator.  }
\end{rema}

{}From Remark \ref{formalinvariance}, for any $\mu\in {\mathbb
C}/{\mathbb Z}$ and logarithmic intertwining operator ${\cal Y}^\mu$
as in (\ref{log:c+n}), the linear map defined by
\[
w_{(1)}\otimes w_{(2)}\mapsto {{\cal Y}^\mu}(w_{(1)},e^{2\pi ip}x)
w_{(2)}= \sum_{\bar n=\mu}\sum_{k\in {\mathbb N}}
{w_{(1)}}_{n;\,k}^{\cal Y}w_{(2)}e^{2\pi ip(-n-1)}x^{-n-1}(\log x+
2\pi ip)^k
\]
is also a logarithmic intertwining operator.  In the strongly graded
case, if the operator in (\ref{log:c+n}) is grading-compatible, then
so is this one.  The right-hand side above can be written as
\begin{eqnarray*}
\lefteqn{e^{-2\pi ip\mu}\sum_{\bar n=\mu}\sum_{k\in {\mathbb
N}}{w_{(1)}}_{n;\,k}^{\cal Y}w_{(2)}x^{-n-1}\sum_{t\in {\mathbb
N}}{k\choose t}(\log x)^{k-t}(2\pi ip)^t}\\
&&=e^{-2\pi ip\mu}\sum_{t\in {\mathbb N}}(2\pi ip)^t\sum_{k\in {\mathbb
N}}{k+t\choose t}\sum_{\bar n=\mu}{w_{(1)}}_{n;\,k+t}^{\cal
Y}w_{(2)}x^{-n-1}(\log x)^k
\end{eqnarray*}
(the coefficient of each power of $x$ being a finite sum over $t$ and
$k$).  We now have:

\begin{propo}
Let $W_1$, $W_2$, $W_3$ be generalized modules for a M\"obius (or
conformal) vertex algebra $V$, and let ${\cal Y}(\cdot, x)\cdot$ be a
logarithmic intertwining operator of type ${W_3\choose W_1\,W_2}$.
For $\mu\in {\mathbb C}/{\mathbb Z}$ and $t \in {\mathbb N}$, define
\[
{\cal X}^{\mu}_t:W_1\otimes W_2\to W_3[\log x]\{x\}
\]
by:
\[
{\cal X}^{\mu}_t: w_{(1)}\otimes w_{(2)}\mapsto \sum_{k\in {\mathbb N}}
{k+t\choose t}\sum_{\bar n=\mu}{w_{(1)}}_{n; \,k+t}^{\cal Y}
w_{(2)}x^{-n-1}(\log x)^k.
\]
Then each ${\cal X}^{\mu}_t$ is a logarithmic intertwining operator of
type ${W_3\choose W_1\,W_2}$.  Equivalently, the operator ${\cal X}_t$
defined by
\begin{equation}\label{newio}
{\cal X}_t: w_{(1)}\otimes w_{(2)}\mapsto \sum_{k\in {\mathbb N}}
{k+t\choose t}\sum_{n\in{\mathbb C}}{w_{(1)}}_{n; \,k+t}^{\cal Y}
w_{(2)}x^{-n-1}(\log x)^k
\end{equation}
is a logarithmic intertwining operator of the same type.  In the
strongly graded case, if ${\cal Y}$ is grading-compatible, then so are
${\cal X}^{\mu}_t$ and ${\cal X}_t$.
\end{propo}
\pf The equivalence of the two statements follows from Remark
\ref{log:mu}.

We now prove that for each $s \in {\mathbb N}$, ${\cal X}^{\mu}_s$ is
a logarithmic intertwining operator.  The conditions (\ref{log:map})
and (\ref{log:ltc}) for ${\cal X}^{\mu}_s$ follow from these
conditions for ${\cal Y}$.  {}From the above we see that for any
integer $p$, $e^{-2\pi ip\mu}\sum_{t\in {\mathbb N}}(2\pi ip)^t{\cal
X}^{\mu}_t$, and hence
\begin{equation}\label{sumx}
{\cal Y}^{\mu}_p = \sum_{t\in {\mathbb N}}(2\pi ip)^t{\cal X}^{\mu}_t,
\end{equation}
is a logarithmic intertwining operator.

To prove (\ref{log:jacobi}), (\ref{log:L(-1)dev}) and
(\ref{log:L(j)b}) for ${\cal X}^{\mu}_s$, we observe that for fixed
$w_{(1)}$, $w_{(2)}$ and $n$ such that $\bar n=\mu$, and in fact for
any finite set of $w_{(1)}$, $w_{(2)}$ and $n$ such that $\bar n=\mu$,
there exists $S \in {\mathbb N}$ with $S \geq s$ such that the
coefficient
\[
{\cal X}^{\mu}_t(w_{(1)} \otimes w_{(2)})_{[-n-1]}
\]
of $x^{-n-1}$ in ${\cal X}^{\mu}_t(w_{(1)} \otimes w_{(2)})$ is $0$ if
$t>S$.  In particular, for any $p \in {\mathbb Z}$, the expression
given by (\ref{sumx}) for the coefficient
\[
{\cal Y}^{\mu}_p(w_{(1)} \otimes w_{(2)})_{[-n-1]}
\]
of $x^{-n-1}$ in ${\cal Y}^{\mu}_p(w_{(1)} \otimes w_{(2)})$ becomes
\[
{\cal Y}^{\mu}_p(w_1 \otimes w_2)_{[-n-1]}
= \sum_{t=0}^{S} (2\pi ip)^t {\cal X}^{\mu}_t(w_{(1)} \otimes w_{(2)})_{[-n-1]}.
\]
Choosing $p \in {\mathbb N}$ and inverting the (finite) Vandermonde
matrix
\[
((2\pi ip)^t)_{p,t=0,\dots,S},
\]
we obtain that for $t \leq S$ (including $t = s$),
\[
{\cal X}^{\mu}_t(w_{(1)} \otimes w_{(2)})_{[-n-1]}
= \sum_{p=0}^{S} a_{tp}^{(S)} {\cal Y}^{\mu}_p(w_{(1)} \otimes w_{(2)})_{[-n-1]}
\]
for some $a_{tp}^{(S)} \in {\mathbb C}$; note that $\sum_{p=0}^{S}
a_{tp}^{(S)} {\cal Y}^{\mu}_p$ is a logarithmic intertwining operator.
(Also note, incidentally, that when $S$ is increased, this formula
remains true even though the numbers $a_{tp}^{(S)}$ change; the
inverses of the (finite) Vandermonde matrices do not stablize as $S$
grows.)  {}From this formula we see easily that the properties
(\ref{log:jacobi}), (\ref{log:L(-1)dev}) and (\ref{log:L(j)b}) for
each ${\cal Y}^{\mu}_p$ imply these properties for ${\cal X}^{\mu}_s$,
and so ${\cal X}^{\mu}_s$ is a logarithmic intertwining operator.  (A
different proof that each ${\cal X}_t$, and hence ${\cal X}^{\mu}_t$,
is a logarithmic intertwining operator is given in the course of the
next three remarks.)

The last assertion is clear.  \epf

\begin{rema}\label{log:fcf}{\rm
Let $W_i$, $W^i$, $i=1,2,3$, be generalized modules for a M\"obius (or
conformal) vertex algebra $V$. If ${\cal Y}(\cdot,x)\cdot$ is a
logarithmic intertwining operator of type ${W_3\choose W_1\,W_2}$ and
\[
\sigma_1: W^1\to W_1,\; \sigma_2: W^2\to W_2 \;\;\mbox{and} \;\; \sigma_3:
W_3\to W^3
\]
are $V$-module homomorphisms, then it is easy to see that
$\sigma_3{\cal Y}(\sigma_1\cdot,x)\sigma_2\cdot$ is a logarithmic
intertwining operator of type ${W^3\choose W^1\,W^2}$.  In the
strongly graded case, if ${\cal Y}$ is grading-compatible, then so is
$\sigma_3{\cal Y}(\sigma_1\cdot,x)\sigma_2\cdot$ (recall {}from Remark
\ref{homsaregradingpreserving} that each $\sigma_j$ preserves the
$\tilde A$-grading).  That is, in categorical language, let ${\cal C}$
be a full subcategory of either the category of $V$-modules, or the
category of generalized $V$-modules, or ${\cal M}_{sg}$ (the category
of strongly graded $V$-modules; recall Notation \ref{MGM}), or ${\cal
GM}_{sg}$ (the category of strongly graded generalized $V$-modules;
again recall Notation \ref{MGM}).  Then the correspondence {}from
${\cal C}\times{\cal C}\times{\cal C}$ to the category ${\bf Vect}$ of
vector spaces given by
\[
(W_1, W_2, W_3)\mapsto {\cal V}^{W_3}_{W_1\,W_2}
\]
is functorial in the third slot and cofunctorial
in the first two slots.  }
\end{rema}

\begin{rema}\label{log:newiorm}{\rm
Now recall {}from Remark \ref{set:L(0)s} that $L(0)-L(0)_s$ commutes
with the actions of both $V$ and ${\mathfrak s}{\mathfrak l}(2)$. So
$(L(0)-L(0)_s)^i$ is a $V$-module homomorphism {}from a generalized
module to itself for any $i\in {\mathbb N}$, and this remains true in
the strongly graded case.  Hence by Remark \ref{log:fcf}, given any
logarithmic intertwining operator ${\cal Y}(\cdot,x)\cdot$ as in
(\ref{log:map0}) and any $i,j,k\in {\mathbb N}$,
\begin{equation}\label{newio'}
(L(0)-L(0)_s)^k{\cal Y}((L(0)-L(0)_s)^i\cdot,x) (L(0)-L(0)_s)^j\cdot
\end{equation}
is again a logarithmic intertwining operator, and in the
strongly graded case, if ${\cal Y}$ is grading-compatible, so is this
operator.  In the next remark we will see that the logarithmic
intertwining operators (\ref{newio}) are just linear combinations of
these.  }
\end{rema}

\begin{rema}{\rm
Let $W_1$, $W_2$, $W_3$ and ${\cal Y}$ be as above and let
$w_{(1)}\in{W_1}_{[n_1]}$ and $w_{(2)}\in{W_2}_{[n_2]}$ for some
complex numbers $n_1$ and $n_2$. Fixing $n\in{\mathbb C}$ and using the notation
$T_{p,k,q}$ in (\ref{log:Tprq}) with $a=n_1$ and $b=n_2$, we rewrite
formula (\ref{log:L(0)^1}) in the proof of Lemma \ref{log:lemma} as
\[
(k+1)T_{p,k+1,q}=(L(0)-n_1-n_2+n+1)T_{p,k,q}-T_{p+1,k,q}-T_{p,k,q+1}
\]
for any $p,q,k\in {\mathbb N}$.  {}From this, by (\ref{log:expand}) we see
that for any $t\in {\mathbb N}$,
\[
{k+t\choose t}T_{p,k+t,q}=\sum_{i,j,l\in {\mathbb N},\;i+j+l=t}
\frac{1}{i!j!l!}(-1)^{i+j}(L(0)-n_1-n_2+n+1)^l T_{p+i,k,q+j}.
\]
Setting $p=q=0$ we get
\begin{eqnarray}\label{log:r+t=?}
\lefteqn{{k+t\choose t}{w_{(1)}}^{\cal Y}_{n;\,k+t}w_{(2)}=
\sum_{i,j,l\in {\mathbb N},\;i+j+l=t}\frac{1}{i!j!l!}(-1)^{i+j}\cdot}\nno\\
&&\cdot(L(0)-n_1-n_2+n+1)^l ((L(0)-n_1)^i{w_{(1)}})^{\cal Y}_{n;\,k}
(L(0)-n_2)^jw_{(2)})\nno\\
\end{eqnarray}
for any $t\in {\mathbb N}$.  (Note that this formula gives an
alternate proof of Proposition \ref{log:logwt}(c).)  Now multiplying
by $x^{-n-1}(\log x)^k$ and then summing over $n\in{\mathbb C}$ and
over $k\in{\mathbb N}$ we see that for every $t\in{\mathbb N}$,
${\cal X}_t$ in (\ref{newio}) is a linear combination of logarithmic
intertwining operators of the form (\ref{newio'}).
}
\end{rema}

In preparation for generalizing basic results {}from \cite{FHL} on
intertwining operators to the logarithmic case, we need to generalize
more of the basic tools.  We now define operators ``$x^{\pm L(0)}$''
for generalized modules, in the natural way:

\begin{defi}{\rm
Let $W$ be a generalized module for a M\"obius (or conformal) vertex
algebra. We define
\[
x^{\pm L(0)}: W\to W\{x\}[\log x]\subset W[\log x]\{x\}
\]
as follows: For any $w\in W_{[n]}$ ($n\in{\mathbb C}$), define
\begin{equation}\label{log:xpmL}
x^{\pm L(0)}w=x^{\pm n}e^{\pm\log x(L(0)-n)}w
\end{equation}
(note that the local nilpotence of $L(0)-n$ on $W_{[n]}$ insures that
the formal exponential series terminates) and then extend linearly to
all $w\in W$.  (Of course, we could also write
\begin{equation}\label{log:xpmLs}
x^{\pm L(0)}=x^{\pm L(0)_s}e^{\pm\log x(L(0)-L(0)_s)},
\end{equation}
using the notation $L(0)_s$.)  We also define operators $x^{\pm
L'(0)}$ on $W^*$ by the condition that for all $w'\in W^*$ and $w\in
W$,
\begin{equation}\label{log:xpmL'}
\langle x^{\pm L'(0)}w',w\rangle=\langle w',x^{\pm L(0)}w\rangle\;
(\in {\mathbb C}\{x\}[\log x]),
\end{equation}
so that
\[
x^{\pm L'(0)}: W^*\to W^*\{x\}[[\log x]].
\]  }
\end{defi}

\begin{rema}\label{3.33}{\rm
Note that these definitions are (of course) compatible with the usual
definitions if $W$ is just an (ordinary) module.  In formula
(\ref{log:xpmL}), $x^{\pm L(0)}$ is defined in a naturally factored
form reminiscent of the factorization invoked in Remark \ref{log:[[]]}
(providing counterexamples there); the symbol $x^{\pm (L(0)-n)}$ is
given meaning by its replacement by $e^{\pm\log x(L(0)-n)}$.  Also
note that in case both $V$ and $W$ are strongly graded, the definition
of $x^{\pm L'(0)}$ given by (\ref{log:xpmL'}), when applied to the
subspace $W'$ of $W^*$, coincides with the definition of $x^{\pm
L'(0)}$ given by (\ref{log:xpmL}) induced {}from the contragredient
module action of $V$ on $W'$. (Recall Theorem \ref{set:W'}.) }
\end{rema}

\begin{rema}{\rm
Note that for $w\in W_{[n]}$, by definition we have
\begin{equation}\label{log:x^L(0)}
x^{\pm L(0)}w=x^{\pm n}\sum_{i\in {\mathbb N}}\frac{(L(0)-n)^iw}{i!}(\pm
\log x)^i\in x^{\pm n}W_{[n]}[\log x].
\end{equation}
It is also handy to have that for any $w\in W$,
\begin{equation}\label{log:inv}
x^{L(0)}x^{-L(0)}w=w=x^{-L(0)}x^{L(0)}w,
\end{equation}
which is clear {}from definition.  Later we will also need the
formula
\begin{equation}\label{log:dx^}
\frac{d}{dx}x^{\pm L(0)}w=\pm x^{-1}x^{\pm L(0)}L(0)w
\end{equation}
for any $w\in W$, i.e.,
\begin{equation}
x\frac{d}{dx}x^{\pm L(0)}w=\pm x^{\pm L(0)}L(0)w,
\end{equation}
or equivalently,
\begin{equation}
\left(x\frac{d}{dx}\mp L(0)\right)x^{\pm L(0)}w=0
\end{equation}
(cf. Lemma \ref{log:de} and Remark \ref{log:[[]]}).  This can be
proved by directly checking that for $w$ homogeneous of generalized
weight $n$,
\[
\frac{d}{dx}e^{\pm\log x(L(0)-n)}w=\pm x^{-1}e^{\pm \log x(L(0)-n)}
(L(0)-n)w,
\]
and hence for such $w$,
\begin{eqnarray*}
\frac{d}{dx}x^{\pm L(0)}w&=&\frac{d}{dx}(x^{\pm n}e^{\pm\log x(L(0)-n)}w)\\
&=&\pm nx^{\pm n-1}e^{\pm\log x(L(0)-n)}w\pm x^{\pm n-1}e^{\pm \log x
(L(0)-n)}(L(0)-n)w\\
&=&\pm x^{\pm n-1}e^{\pm\log x(L(0)-n)}L(0)w=\pm x^{-1}x^{\pm L(0)}L(0)w.
\end{eqnarray*}
}
\end{rema}

In the statement (and proof) of the next result, we shall use
expressions of the type
\begin{eqnarray*}
(1-x)^{L(0)}&=&\sum_{k\in \mathbb{N}}{L(0)\choose k}(-x)^{k}\nn
&=&\sum_{k\in \mathbb{N}}\frac{L(0)(L(0)-1)\cdots
(L(0)-k+1)}{k!}(-x)^{k},
\end{eqnarray*}
which also equals
\[
e^{L(0)\log (1-x)},
\]
as well as expressions involving
\[
(x(1-yx)^{-1})^{n}=\sum_{k\in \mathbb{N}}{-n\choose k}x^{n}(-yx)^{k}
\]
for $n \in \C$ and
\begin{eqnarray*}
\log (x(1-yx)^{-1})&=&\log x+\log (1-yx)^{-1}\nn
&=&\log x+\sum_{k\ge 1}\frac{1}{k}(yx)^{k}.
\end{eqnarray*}

We can now state and prove generalizations to logarithmic intertwining
operators of three standard formulas for (ordinary) intertwining
operators, namely, formulas (5.4.21), (5.4.22) and (5.4.23) of
\cite{FHL}.  The result is (see also \cite{Mi} for Parts (a) and (b)):

\begin{propo}
Let ${\cal Y}$ be a logarithmic intertwining operator of type
${W_3\choose W_1\,W_2}$ and let $w\in W_1$. Then
\begin{description}
\item{(a)}
\begin{equation}\label{log:p1}
e^{yL(-1)}{\cal Y}(w,x)e^{-yL(-1)}={\cal Y}(e^{yL(-1)}w,x)={\cal Y}(w,x+y)
\end{equation}
(recall (\ref{log:not1}))
\item{(b)}
\begin{equation}\label{log:p2}
y^{L(0)}{\cal Y}(w,x)y^{-L(0)}={\cal Y}(y^{L(0)}w,xy)
\end{equation}
(recall (\ref{log:not3}))
\item{(c)}
\begin{equation}\label{log:p3}
e^{yL(1)}{\cal Y}(w,x)e^{-yL(1)}={\cal
Y}(e^{y(1-yx)L(1)}(1-yx)^{-2L(0)}w,x(1-yx)^{-1}).
\end{equation}
\end{description}
\end{propo}
\pf {}From (\ref{log:L(j)b}) with $j=-1$ we see that for any $w\in W_1$,
\[
L(-1){\cal Y}(w,x)={\cal Y}(L(-1)w,x)+{\cal Y}(w,x)L(-1).
\]
This implies
\[
\frac{y^n(L(-1))^n}{n!}{\cal Y}(w,x)=\sum_{i,j\in{\mathbb N},\;i+j=n}{\cal
Y} \bigg(\frac {y^i(L(-1))^i}{i!}w,x\bigg)\frac{y^j(L(-1))^j}{j!}
\]
for any $n\in{\mathbb N}$, where $y$ is a new formal
variable. Summing over $n\in{\mathbb N}$ we see that for any $w\in W_1$,
\[
e^{yL(-1)}{\cal Y}(w,x)={\cal Y}(e^{yL(-1)}w,x)e^{yL(-1)},
\]
and hence by (\ref{log:ck1}),
\begin{equation}\label{log:p1-1}
e^{yL(-1)}{\cal Y}(w,x)e^{-yL(-1)}={\cal Y}(e^{yL(-1)}w,x)
={\cal Y}(w,x+y)
\end{equation}
(and this expression also equals $e^{y\frac d{dx}}{\cal Y}(w,x)$).
Note that all the expressions in (\ref{log:p1-1}) remain well defined
if we replace $y$ by any element of $y\C[x][[y]]$ and that
(\ref{log:p1-1}) still holds if we replace $y$ by any such element.
(But note that (\ref{log:ck1}) would fail to hold if we replaced $y$
by for example $yx \in y\C[x][[y]]$.)

For (b), note that for homogeneous $w_{(1)}\in {W_1}$ and $w_{(2)}\in
{W_2}$, by (\ref{log:e^L(0)}) with the formal variable $y$ replaced by
the formal variable $\log y$, we get (recalling Proposition
\ref{log:logwt}(b) and (\ref{log:xpmL}))
\begin{eqnarray*}
\lefteqn{y^{L(0)}({w_{(1)}}_{n;\,k}^{\cal Y}w_{(2)})=}\\
&&\sum_{l\in{\mathbb N}}{k+l\choose k}(y^{L(0)}w_{(1)})_{n;\,k+l}^{\cal
Y}(y^{L(0)}w_{(2)})y^{-n-1}(\log y)^l
\end{eqnarray*}
for any $n\in{\mathbb C}$ and $k\in{\mathbb N}$.
Multiplying this by $x^{-n-1}(\log x)^k$, summing over $n\in
{\mathbb C}$ and $k\in {\mathbb N}$ and using (\ref{log:not3}) we get
\[
y^{L(0)}{\cal Y}(w_{(1)},x)w_{(2)}
={\cal Y}(y^{L(0)}w_{(1)},xy)y^{L(0)}w_{(2)}.
\]
Formula (\ref{log:p2}) then follows {}from (\ref{log:inv}).

Finally we prove (c).  {}From (\ref{log:L(j)b}) with $j=1$ we see that
for any $w\in W_1$,
\[
L(1){\cal Y}(w,x)={\cal Y}((L(1)+2xL(0)+x^2L(-1))w,x)+
{\cal Y}(w,x)L(1).
\]
This implies that
\[
e^{yL(1)}{\cal Y}(w,x)={\cal Y}(e^{y(L(1)+2xL(0)+x^2L(-1))}w,x)
e^{yL(1)},
\]
or
\[
e^{yL(1)}{\cal Y}(w,x)e^{-yL(1)}={\cal Y}
(e^{y(L(1)+2xL(0)+x^2L(-1))}w,x).
\]
Using the identity
\[
e^{y(L(1)+2xL(0)+x^2L(-1))}=e^{yx^2(1-yx)^{-1}L(-1)}
e^{y(1-yx)L(1)}(1-yx)^{-2L(0)},
\]
whose proof is exactly the same as the proof of formula (5.2.41) of
\cite{FHL}, we obtain
\begin{equation}\label{log:p4}
e^{yL(1)}{\cal Y}(w,x)e^{-yL(1)}={\cal
Y}(e^{yx^2(1-yx)^{-1}L(-1)}e^{y(1-yx)L(1)}(1-yx)^{-2L(0)}w,x).
\end{equation}
But by (\ref{log:p1-1}) with $y$ replaced by $yx^2(1-yx)^{-1}$, 
the right-hand side of (\ref{log:p4}) 
is equal to 
\begin{equation}\label{log:p4-r}
{\cal Y}(e^{y(1-yx)L(1)}(1-yx)^{-2L(0)}w,x+yx^2(1-yx)^{-1}).
\end{equation}
Since 
\[
(x+yx^2(1-yx)^{-1})^{n}=(x(1-yx)^{-1})^{n}
\]
for $n \in \C$ and
\[
\log (x+yx^2(1-yx)^{-1})=\log (x(1-yx)^{-1}),
\]
(\ref{log:p4-r}) is equal to the right-hand side of 
(\ref{log:p3}), proving (\ref{log:p3}). 
\epf

\begin{rema}{\rm
The following formula, also a generalization of the corresponding
formula in the ordinary case (see (\ref{xL(0)L(j)})), will be needed:
For $j=-1, 0, 1$,
\begin{equation}\label{log:xLx^}
x^{L(0)}L(j)x^{-L(0)}=x^{-j}L(j).
\end{equation}
To prove this, we first observe that for any $m\in{\mathbb C}$,
$[L(0)-m,L(j)]=-jL(j)$ implies that
\[
e^{\log x(L(0)-m)}L(j)e^{-\log x(L(0)-m)}=e^{-j\log x}L(j).
\]
Hence, for a generalized module element $w$ homogeneous of generalized
weight $n$,
\begin{eqnarray*}
x^{L(0)}L(j)w&=&x^{n-j}e^{\log x(L(0)-n+j)}L(j)w\\
&=&x^{n-j}e^{-j\log x}L(j)e^{\log x(L(0)-n+j)}w\\
&=&x^{n-j}L(j)e^{\log x(L(0)-n)}w\\
&=&x^{-j}L(j)x^{L(0)}w,
\end{eqnarray*}
and (\ref{log:xLx^}) then follows immediately {}from (\ref{log:inv}).
(Another proof: By Remark \ref{set:L(0)s}, the second factor on the
right-hand side of (\ref{log:xpmLs}) commutes with $L(j)$.)
}
\end{rema}

\begin{rema}{\rm {}From (\ref{log:xLx^}) we see that
\begin{equation}\label{xe^Lx}
x^{L(0)}e^{yL(j)}x^{-L(0)}=e^{yx^{-j}L(j)}.
\end{equation}
}
\end{rema}

For an ordinary module $W$ for a vertex operator algebra and any
$a\in{\mathbb C}$, the operator $e^{aL(0)}$ on $W$ is defined
by
\begin{equation}\label{eaL0ordinary}
e^{aL(0)}w=e^{ah}w
\end{equation}
for any homogeneous $w\in W_{(h)}$, $h\in{\mathbb C}$ and then by linear
extension to any $w\in W$.  More generally, for a generalized module
$W$ for a M\"obius (or conformal) vertex algebra and any $a\in{\mathbb
C}$, we define the operator $e^{aL(0)}$ on $W$ by
\begin{equation}\label{eaL0}
e^{aL(0)}w=e^{ah}e^{a(L(0)-h)}w
\end{equation}
for any homogeneous $w\in W_{[h]}$, $h\in{\mathbb C}$ and then by
linear extension to all $w\in W$. (Note that for a formal variable
$x$, we already have $e^{xL(0)}w=e^{hx}e^{x(L(0)-h)}w$.)  {}From the
definition,
\begin{equation}
e^{aL(0)}e^{-aL(0)}w=w.
\end{equation}
Recalling Remark
\ref{set:L(0)s} for the notation $L(0)_s$, we see that
\begin{equation}\label{eaL0-general}
e^{aL(0)}=e^{aL(0)_s}e^{a(L(0)-L(0)_s)}\;\;\mbox{on}\;W,
\end{equation}
where the exponential series $e^{a(L(0)-L(0)_s)}$ terminates on each
element of $W$.

\begin{rema}\label{analyticallyconvergent}{\rm
The operators defined in (\ref{eaL0ordinary}) and (\ref{eaL0}) can be
alternatively defined or viewed as the (analytically) convergent sums
of the indicated exponential series of operators; these operators act
on the (finite-dimensional) subspaces of $W$ generated by the repeated
action of $L(0)$ on homogeneous vectors $w \in W$.
}
\end{rema}

\begin{rema}\label{exponentialaVhom}{\rm
The operator $e^{a(L(0)-L(0)_s)}$ on $W$ is a $V$-homomorphism, in
view of Remark \ref{set:L(0)s} (cf.\ Remark \ref{log:newiorm}).  Let
$r$ be an integer.  Then $e^{2\pi irL(0)_s}$ is also a
$V$-homomorphism, by Remark \ref{congruent}.  Thus for $r \in \Z$,
$e^{2\pi irL(0)}$ is a $V$-homomorphism.  In the strongly graded
case, all of these $V$-homomorphisms are grading-preserving.}
\end{rema}

\begin{rema}{\rm
We now recall some identities about the action of ${\mathfrak s}{\mathfrak l}
(2)$ on any of its modules. For convenience we put them in the
following form:
\begin{equation}\label{log:SL2-1}
e^{xL(-1)}\left(\begin{array}{c}L(-1)\\L(0)\\L(1)\end{array}\right)
e^{-xL(-1)}=\left(\begin{array}{ccc}1&0&0\\-x&1&0\\x^2&-2x&1\end{array}
\right)\left(\begin{array}{c}L(-1)\\L(0)\\L(1)\end{array}\right)
\end{equation}
\begin{equation}\label{log:SL2-2}
e^{xL(0)}\left(\begin{array}{c}L(-1)\\L(0)\\L(1)\end{array}\right)
e^{-xL(0)}=\left(\begin{array}{ccc}e^x&0&0\\0&1&0\\0&0&e^{-x}\end{array}
\right)\left(\begin{array}{c}L(-1)\\L(0)\\L(1)\end{array}\right)
\end{equation}
\begin{equation}\label{log:SL2-3}
e^{xL(1)}\left(\begin{array}{c}L(-1)\\L(0)\\L(1)\end{array}\right)
e^{-xL(1)}=\left(\begin{array}{ccc}1&2x&x^2\\0&1&x\\0&0&1\end{array}
\right)\left(\begin{array}{c}L(-1)\\L(0)\\L(1)\end{array}\right).
\end{equation}
Formula (\ref{log:SL2-2}) follows {}from (5.2.12) and (5.2.13) in
\cite{FHL}; Formula (\ref{log:SL2-3}) follows {}from (5.2.14) in
\cite{FHL} and $[L(1),L(-1)]=2L(0)$; and formula (\ref{log:SL2-1})
follows {}from (\ref{log:SL2-3}) and the fact that
\[
L(-1)\mapsto L(1),\; L(0)\mapsto-L(0),\; L(1)\mapsto L(-1)
\]
is a Lie algebra automorphism of ${\mathfrak s}{\mathfrak l}(2)$.
}
\end{rema}

\begin{rema}\label{log:Lj2rema}{\rm
It is convenient to note that the ${\mathfrak s}{\mathfrak l}(2)$-bracket
relations (\ref{log:L(j)b}) are equivalent to
\begin{equation}\label{log:L(j)b2}
{\cal Y}(L(j)w_{(1)},x)=\sum_{i=0}^{j+1}{j+1\choose i}(-x)^i
[L(j-i),{\cal Y}(w_{(1)},x)]
\end{equation}
for $j=-1$, $0$ and $1$. This can be easily checked by writing
(\ref{log:L(j)b}) as
\[
\left(\begin{array}{c}[L(-1),{\cal Y}(w_{(1)},x)]\\{}[L(0),{\cal Y}
(w_{(1)},x)]\\{}[L(1),{\cal Y}(w_{(1)},x)]\end{array}\right)
=\left(\begin{array}{ccc}1&0&0\\x&1&0\\x^2&2x&1\end{array}\right)
\left(\begin{array}{c}{\cal Y}(L(-1)w_{(1)},x)\\{\cal Y}
(L(0)w_{(1)},x)\\{\cal Y}(L(1)w_{(1)},x)\end{array}\right)
\]
and then multiplying this by the inverse of the invertible matrix on
the right-hand side, obtained by replacing $x$ by $-x$.  (Of course,
in the case where $V$ is conformal, this equivalence is already
encoded in the symmetry of the Jacobi identity.)  }
\end{rema}

We have already defined the natural process of multiplying the formal
variable in a logarithmic intertwining operator by $e^{2\pi ip}$ for
$p \in \Z$ (recall (\ref{substitutionofe2piipx})), and this process
yields another logarithmic intertwining operator (recall Remark
\ref{formalinvariance}).  It is natural to generalize this
substitution process to that of multiplying the formal variable in a
logarithmic intertwining operator by the exponential $e^{\zeta}$ of
any complex number $\zeta$.  As in the special case $\zeta = 2\pi ip$,
the process will depend on $\zeta$, not just on $e^{\zeta}$, but we
will still find it convenient to use the shorthand symbol $e^{\zeta}$
in our notation for the process.  In Section 7 of \cite{tensor2}, we
introduced this procedure in the case of ordinary (nonlogarithmic)
intertwining operators, and we now carry it out in the general
logarithmic case.  We are about to use this substitution mostly for
\[
\zeta = (2r+1)\pi i,\; r \in \Z.
\]

Let $(W_1,Y_1)$, $(W_2,Y_2)$ and $(W_3,Y_3)$ be generalized modules
for a M\"obius (or conformal) vertex algebra $V$.  Let ${\cal Y}$ be a
logarithmic intertwining operator of type ${W_3\choose W_1\, W_2}$.
For any complex number $\zeta$ and any $w_{(1)}\in W_1$ and
$w_{(2)}\in W_2$, set
\begin{equation}\label{log:subs}
{\cal Y} (w_{(1)},e^\zeta x)w_{(2)} = {\cal
Y}(w_{(1)},y)w_{(2)}\lbar_{y^n=e^{\zeta n}x^n,\; (\log y)^k=(\zeta+\log
x)^k,\;n\in{\mathbb C},\;k\in{\mathbb N}},
\end{equation}
a well-defined element of $W_3[\log x]\{x\}$.  Note that this element
indeed depends on $\zeta$, not just on $e^\zeta$.

\begin{rema}{\rm
In Section 4 below we will take the further step of specializing the
formal variable $x$ to $1$ (or equivalently, the formal variable $y$
to $e^{\zeta}$) in (\ref{log:subs}); that is, we will consider ${\cal
Y} (w_{(1)},e^\zeta)w_{(2)}$.
}
\end{rema}

Given any $r\in {\mathbb Z}$, we define
\[
\Omega_r({\cal Y}): W_2\otimes W_1\to W_3[\log x]\{x\}
\]
by the formula
\begin{equation}\label{Omega_r}
\Omega_{r}({\cal Y})(w_{(2)},x)w_{(1)} = e^{xL(-1)} {\cal
Y}(w_{(1)},e^{ (2r+1)\pi i}x)w_{(2)}
\end{equation}
for $w_{(1)}\in W_{1}$ and $w_{(2)}\in W_{2}$.  This expression is
indeed well defined because of the truncation condition
(\ref{log:ltc}) (recall Remark \ref{ordinaryandlogintwops}).  The
following result generalizes Proposition 7.1 of \cite{tensor2} (for
the ordinary intertwining operator case) and has essentially the same
proof as that proposition, which in turn generalized Proposition 5.4.7
of \cite{FHL}:

\begin{propo}\label{log:omega}
The operator $\Omega_r({\cal Y})$ is a logarithmic intertwining
operator of type ${W_3\choose W_2\, W_1}$. Moreover,
\begin{equation}\label{log:or}
\Omega_{-r-1}(\Omega_r({\cal Y}))=\Omega_r(\Omega_{-r-1}({\cal Y}))
={\cal Y}.
\end{equation}
In the strongly graded case, if ${\cal Y}$ is grading-compatible, then
so is $\Omega_r({\cal Y})$, and in particular, the correspondence
${\cal Y}\mapsto \Omega_r({\cal Y})$ defines a linear isomorphism {}from
${\cal V}_{W_1\,W_2}^{W_3}$ to ${\cal V}_{W_2\,W_1}^{W_3}$, and we
have
\[
 N_{W_1\,W_2}^{W_3}=N_{W_2\,W_1}^{W_3}.
\]
\end{propo}
\pf The lower truncation condition (\ref{log:ltc}) is clear. {From}
the Jacobi identity (\ref{log:jacobi}) for ${\cal Y}$,
\begin{eqnarray}
\lefteqn{\dps x^{-1}_0\delta \left( {x_1-y\over x_0}\right)
Y_3(v,x_1){\cal Y}(w_{(1)},y)w_{(2)}}\nno\\
&&\hspace{2em}- x^{-1}_0\delta \left( {y-x_1\over -x_0}\right)
{\cal Y}(w_{(1)},y)Y_2(v,x_1)w_{(2)}\nno \\
&&{\dps = y^{-1}\delta \left( {x_1-x_0\over y}\right)
{\cal Y}(Y_1(v,x_0)w_{(1)},y)
w_{(2)}},
\end{eqnarray}
with $v \in V$, $w_{(1)} \in W_1$ and $w_{(2)} \in W_2$, we obtain
\begin{eqnarray}\label{75}
\dps x^{-1}_0\delta \left( {x_1-y\over x_0}\right)
e^{-yL(-1)}Y_3(v,x_1){\cal Y}(w_{(1)},y)w_{(2)}
\lbar_{y^n=e^{(2r+1)\pi in}x_2^n,\; (\log y)^k=((2r+1)\pi i+\log
x_2)^k,\;n\in{\mathbb C},\;k\in{\mathbb N}}
\hspace{1.5em}\nno\\
- x^{-1}_0\delta \left( {y-x_1\over -x_0}\right)
e^{-yL(-1)}{\cal Y}(w_{(1)},y)Y_2(v,x_1)w_{(2)}
\lbar_{y^n=e^{(2r+1)\pi in}x_2^n,\; (\log y)^k=((2r+1)\pi i+\log
x_2)^k,\;n\in{\mathbb C},\;k\in{\mathbb N}}\nno \\
{\dps = y^{-1}\delta \left( {x_1-x_0\over y}\right)
e^{-yL(-1)}{\cal Y}(Y_1(v,x_0)w_{(1)},y)
w_{(2)}\lbar_{y^n=e^{(2r+1)\pi in}x_2^n,\; (\log y)^k=((2r+1)\pi i+\log
x_2)^k,\;n\in{\mathbb C},\;k\in{\mathbb N}}}.
\end{eqnarray}
The first term of the left-hand side of (\ref{75}) is equal to
\begin{eqnarray*}
&{\dps x^{-1}_0\delta \left( {x_1-y\over x_0}\right)
Y_3(v,x_1-y)e^{-yL(-1)}{\cal Y}(w_{(1)},y)w_{(2)}
\lbar_{y^n=e^{(2r+1)\pi in}x_2^n,\; (\log y)^k=((2r+1)\pi i+\log
x_2)^k,\;n\in{\mathbb C},\;k\in{\mathbb N}}}&\nno\\
&{\dps =x^{-1}_0\delta \left( {x_1+x_{2}\over x_0}\right)
Y_3(v,x_0)\Omega_{r}({\cal Y})(w_{(2)},x_{2})w_{(1)},}&
\end{eqnarray*}
the second term  is equal to
$$
- x^{-1}_0\delta \left( {-x_{2}-x_1\over -x_0}\right)
\Omega_{r}({\cal Y})(Y_2(v,x_1)w_{(2)},x_{2})w_{(1)}
$$
and the right-hand side of (\ref{75}) is equal to
$$
-x_{2}^{-1}\delta \left( {x_1-x_0\over -x_{2}}\right)
\Omega_{r}({\cal Y})(w_{(2)},x_{2})Y_1(v,x_0)w_{(1)}.
$$
Substituting  into (\ref{75}) we obtain
\begin{eqnarray}
\lefteqn{x^{-1}_0\delta \left( {x_1+x_{2}\over x_0}\right)
Y_3(v,x_0)\Omega_{r}({\cal Y})(w_{(2)},x_{2})w_{(1)}}\nno\\
&&\hspace{2em}- x^{-1}_0\delta \left( {x_{2}+x_1\over x_0}\right)
\Omega_{r}({\cal Y})(Y_2(v,x_1)w_{(2)},x_{2})w_{(1)}\nno\\
&&=-x_{2}^{-1}\delta \left( {x_1-x_0\over -x_{2}}\right)
\Omega_{r}({\cal Y})(w_{(2)},x_{2})Y_1(v,x_0)w_{(1)},
\end{eqnarray}
which is equivalent to
\begin{eqnarray}
\lefteqn{x^{-1}_1\delta \left( {x_0-x_{2}\over x_1}\right)
Y_3(v,x_0)\Omega_{r}({\cal Y})(w_{(2)},x_{2})w_{(1)}}\nno\\
&&\hspace{2em}-x_{1}^{-1}\delta \left( {x_2-x_0\over -x_{1}}\right)
\Omega_{r}({\cal Y})(w_{(2)},x_{2})Y_1(v,x_0)w_{(1)}\nno\\
&&= x^{-1}_2\delta \left( {x_{0}-x_1\over x_2}\right)
\Omega_{r}({\cal Y})(Y_2(v,x_1)w_{(2)},x_{2})w_{(1)}
\end{eqnarray}
(recall (\ref{2termdeltarelation})).  This in turn is the Jacobi
identity for $\Omega_{r}({\cal Y})$ (with the roles of $x_0$ and $x_1$
reversed in (\ref{log:jacobi})).

To prove the $L(-1)$-derivative property (\ref{log:L(-1)dev}) for
$\Omega_{r}({\cal Y})$, first note that {}from (\ref{log:subs}) and the
$L(-1)$-derivative property for ${\cal Y}$,
\begin{eqnarray}
\frac{d}{dx}{\cal Y} (w_{(1)},e^{\zeta} x)w_{(2)}
&=& e^{\zeta}\left(\frac{d}{dy} {\cal
Y}(w_{(1)},y)w_{(2)}\right)\lbar_{y^n=e^{\zeta n}x^n,\; (\log
y)^k=(\zeta+\log x)^k,\;n\in{\mathbb C},\;k\in{\mathbb N}}\nno\\
&=& e^{\zeta}{\cal Y} (L(-1)w_{(1)},e^\zeta x)w_{(2)},
\end{eqnarray}
and in particular,
\begin{equation}
\frac{d}{dx}{\cal Y} (w_{(1)},e^{(2r+1)\pi i} x)w_{(2)}
= -{\cal Y} (L(-1)w_{(1)},e^{(2r+1)\pi i} x)w_{(2)}.
\end{equation}
Thus, by using formula (\ref{log:L(j)b}) with $j=-1$ we have
\begin{eqnarray}
\lefteqn{\frac{d}{dx}\Omega_{r}({\cal Y})(w_{(2)},x)w_{(1)}=
\frac{d}{dx}e^{xL(-1)}{\cal Y}(w_{(1)},e^{ (2r+1)\pi i}x)w_{(2)}}\nno\\
&&=e^{xL(-1)}L(-1){\cal Y}(w_{(1)},e^{ (2r+1)\pi i}x)w_{(2)}+
e^{xL(-1)}\frac{d}{dx}{\cal Y}(w_{(1)},e^{ (2r+1)\pi i}x)w_{(2)}\nno\\
&&=e^{xL(-1)}L(-1){\cal Y}(w_{(1)},e^{ (2r+1)\pi i}x)w_{(2)}-
e^{xL(-1)}{\cal Y}(L(-1)w_{(1)},e^{ (2r+1)\pi i}x)w_{(2)}\nno\\
&&=e^{xL(-1)}{\cal Y}(w_{(1)},e^{ (2r+1)\pi i}x)L(-1)w_{(2)}\nno\\
&&=\Omega_{r}({\cal Y})(L(-1)w_{(2)},x)w_{(1)},
\end{eqnarray}
as desired.

In the case that $V$ is M\"obius, we prove the ${\mathfrak
s}{\mathfrak l}(2)$-bracket relations (\ref{log:L(j)b}) for
$\Omega_r({\cal Y})$. By using the ${\mathfrak s}{\mathfrak
l}(2)$-bracket relations for ${\cal Y}$ and the relations
\[
e^{xL(-1)}L(j)e^{-xL(-1)}=\sum_{i=0}^{j+1}{j+1\choose i}(-x)^iL(j-i)
\]
for $j=-1, 0$ and $-1$ (see (\ref{log:SL2-1})), we have
\begin{eqnarray*}
\lefteqn{\Omega_r({\cal Y})(L(j)w_{(2)},x)w_{(1)}=
e^{xL(-1)}{\cal Y}(w_{(1)},e^{(2r+1)\pi i}x)L(j)w_{(2)}}\\
&&=e^{xL(-1)}L(j){\cal Y}(w_{(1)},e^{(2r+1)\pi i}x)w_{(2)}\\
&&\hspace{4em}-e^{xL(-1)}\sum_{i=0}^{j+1}{j+1\choose i}(-x)^i{\cal Y}
(L(j-i)w_{(1)},e^{(2r+1)\pi i}x)w_{(2)}\\
&&=\sum_{i=0}^{j+1}{j+1\choose i}(-x)^iL(j-i)e^{xL(-1)}{\cal
Y}(w_{(1)},e^{(2r+1)\pi i}x)w_{(2)}\\
&&\hspace{4em}-e^{xL(-1)}\sum_{i=0}^{j+1}{j+1\choose i}(-x)^i{\cal Y}
(L(j-i)w_{(1)},e^{(2r+1)\pi i}x)w_{(2)}\\
&&=\sum_{i=0}^{j+1}{j+1\choose i}(-x)^i(L(j-i)\Omega_r({\cal
Y})(w_{(2)},x)-\Omega_r({\cal Y})(w_{(2)},x)L(j-i))w_{(1)},
\end{eqnarray*}
which is the alternative form (\ref{log:L(j)b2}) of the ${\mathfrak
s}{\mathfrak l}(2)$-bracket relations for $\Omega_r({\cal Y})$.

The identity (\ref{log:or}) is clear {from} the definitions of
$\Omega_{r}({\cal Y})$ and $\Omega_{-r-1}({\cal Y})$, and the
remaining assertions are clear. \epf

\begin{rema}\label{Ys1s2s3}{\rm
For each triple $s_1, s_2, s_3\in{\mathbb Z}$, the
logarithmic intertwining operator ${\cal Y}$ gives rise to a
logarithmic intertwining
operator ${\cal Y}_{[s_1,s_2, s_3]}$ of the same type, defined by
\[
{\cal Y}_{[s_1,s_2, s_3]}(w_{(1)},x)=e^{2\pi is_1 L(0)}{\cal
Y}(e^{2\pi is_2 L(0)} w_{(1)},x)e^{2\pi is_3 L(0)}
\]
for $w_{(1)}\in W_1$, by Remarks \ref{log:fcf} and
\ref{exponentialaVhom}.  In the strongly graded case, if ${\cal Y}$ is
grading-compatible, so is ${\cal Y}_{[s_1,s_2, s_3]}$.  Clearly,
\[
{\cal Y}_{[0,0,0]}={\cal Y}
\]
and for $r_1,r_2,r_3,s_1,s_2,s_3\in{\mathbb Z}$,
\[
({\cal Y}_{[r_1,r_2, r_3]})_{[s_1,s_2, s_3]}={\cal Y}_{[r_1+s_1,
r_2+s_2, r_3+s_3]}.
\]
For any $a\in{\mathbb C}$, we have the formula
\begin{equation}\label{710}
e^{aL(0)}{\cal Y}(w_{(1)},x)e^{-aL(0)}={\cal Y}
(e^{aL(0)}w_{(1)},e^ax)
\end{equation}
(cf.\ (\ref{log:p2})).  This is proved by imitating the proof of
(\ref{log:p2}), replacing $y^{L(0)}$ by $e^{aL(0)}$, $y$ by $e^a$ and
$\log y$ by $a$ in that proof, using (\ref{eaL0}) in place of
(\ref{log:xpmL}) and keeping in mind formula (\ref{log:subs}).  (When
(\ref{log:e^L(0)}) is used in this proof, for homogeneous elements
$w_{(1)}$ and $w_{(2)}$, the exponential series all terminate, as does
the sum over $l \in {\mathbb{N}}$.)  {}From this, we see that
(\ref{log:or}) generalizes to
\[
\Omega_s(\Omega_r({\cal Y}))={\cal Y}_{[r+s+1,-(r+s+1),-(r+s+1)]}
={\cal Y}(\cdot,e^{2\pi i (r+s+1)}\cdot)
\]
for all $r,s \in \Z$.
}
\end{rema}

In case $V$, $W_1$, $W_2$ and $W_3$ are strongly graded, which we now
assume, we have the concept of ``$r$-contragredient operator'' as
follows (in the ordinary intertwining operator case this was
introduced in \cite{tensor2}): Given a grading-compatible logarithmic
intertwining operator ${\cal Y}$ of type ${W_3\choose W_1\,W_2}$ and
an integer $r$, we define the {\em $r$-contragredient operator of
${\cal Y}$} to be the linear map
\begin{eqnarray*}
W_1\otimes W'_3&\to&W'_2\{x\}[[\log x]]\\
w_{(1)}\otimes w'_{(3)}&\mapsto&A_r({\cal Y})(w_{(1)},x)w'_{(3)}
\end{eqnarray*}
given by
\begin{eqnarray}\label{log:Ardef}
\lefteqn{\langle A_r({\cal Y})(w_{(1)},x)w'_{(3)},w_{(2)}
\rangle_{W_2}}\nno\\
&&=\langle w'_{(3)},{\cal Y}(e^{xL(1)}e^{(2r+1)\pi iL(0)}(x^{-L(0)})^2
w_{(1)},x^{-1})w_{(2)}\rangle_{W_3},
\end{eqnarray}
for any $w_{(1)}\in W_1$, $w_{(2)}\in W_2$ and $w'_{(3)}\in W'_3$,
where we use the notation
\[
f(x^{-1})=\sum_{m\in{\mathbb N},\,n\in{\mathbb C}}w_{n,m} x^{-n}(-\log x)^m
\]
for any
\[
f(x)=\sum_{m\in{\mathbb N},\,n\in {\mathbb C}}w_{n,m}
x^n(\log x)^m\in{\cal W}\{x\}[[\log x]],
\]
${\cal W}$ any vector space
(not involving $x$).  Note that for the case $W_1=V$ and $W_2=W_3=W$,
the operator $A_r({\cal Y})$ agrees with the contragredient vertex
operator ${\cal Y}'$ (recall (\ref{yo}) and (\ref{y'})) for any
$r\in{\mathbb Z}$.

We have the following result generalizing Proposition 7.3 in
\cite{tensor2} for ordinary intertwining operators, and having
essentially the same proof (and also generalizing Theorem 5.5.1 and
Proposition 5.5.2 of \cite{FHL}):

\begin{propo}\label{log:A}
The $r$-contragredient operator $A_r({\cal Y})$ of a
grading-compatible logarithmic intertwining operator ${\cal Y}$ of
type ${W_3\choose W_1\,W_2}$ is a grading-compatible logarithmic
intertwining operator of type ${W'_2\choose W_1\,W'_3}$. Moreover,
\begin{equation}\label{log:ar}
A_{-r-1}(A_r({\cal Y}))=A_r(A_{-r-1}({\cal Y}))={\cal Y}.
\end{equation}
In particular, the correspondence ${\cal Y}\mapsto A_r({\cal Y})$
defines a linear isomorphism {}from ${\cal V}_{W_1\,W_2}^{W_3}$ to
${\cal V}_{W_1\,W'_3}^{W'_2}$, and we have
\[
N_{W_1\,W_2}^{W_3}=N_{W_1\,W'_3}^{W'_2}.
\]
\end{propo}
\pf First we need to show that for $w_{(1)}\in W_1$ and $w'_{(3)}\in
W'_3$,
\begin{equation}\label{finitelymanypowersoflogx}
A_r({\cal Y})(w_{(1)},x)w'_{(3)} \in W'_2[\log x]\{x\},
\end{equation}
that is, for each power of $x$ there are only finitely many powers of
$\log x$.  This and the lower truncation condition (\ref{log:ltc}) as
well as the grading-compatibility condition (\ref{gradingcompatcondn})
follow {}from a variant of the argument proving the lower truncation
condition for contragredient vertex operators (recall
(\ref{truncationforY'})):

Fix elements $w_{(1)}\in W_1$ and $w'_{(3)}\in W'_3$ homogeneous with
respect to the double gradings of $W_1$ and $W'_3$ (it will suffice to
prove the desired assertions for such elements), and in fact take
\[
w_{(1)}\in W_1^{(\beta)} \;\;\mbox{and}\;\; w'_{(3)}\in
(W'_3)^{(\gamma)},
\]
where $\beta$ and $\gamma$ are elements of the abelian group $\tilde
A$, in the notation of Definition \ref{def:dgw}, and fix $n \in \C$.
The right-hand side of (\ref{log:Ardef}) is a (finite) sum of terms of
the form
\begin{equation}\label{w3Yw2}
\langle w'_{(3)},{\cal Y}(w,x^{-1})w_{(2)}\rangle x^p (\log x)^q
\end{equation}
where $w \in W_1$ is doubly homogeneous and in fact $w \in
W_1^{(\beta)}$ (by (\ref{m-L(n)-A})), and where $p \in \C$ and $q \in
{\mathbb N}$.  (The pairing $\langle \cdot,\cdot \rangle$ is between
$W'_3$ and $W_3$.)  Let $w_{(2)}^{(-\beta-\gamma)}$ be the component
of (the arbitrary element) $w_{(2)}$ in $W_2^{(-\beta-\gamma)}$, with
respect to the $\tilde A$-grading.  Then (\ref{w3Yw2}) equals
\begin{equation}\label{w3Yw2betagamma}
\langle w'_{(3)},{\cal Y}(w,x^{-1})w_{(2)}^{(-\beta-\gamma)}\rangle x^p
(\log x)^q
\end{equation}
because of the grading-compatibility condition
(\ref{gradingcompatcondn}) together with (\ref{W'beta}).  (This is why
we need our logarithmic intertwining operators to be
grading-compatible.)  This shows in particular that
\begin{equation}
{w_{(1)}}_{N;\,K}^{A_r({\cal Y})}w'_{(3)} \in (W'_{2})^{(\beta + \gamma)}
\end{equation}
for $N \in \C$ and $K \in {\mathbb N}$, so that
(\ref{gradingcompatcondn}) holds for $A_r({\cal Y})$.  Let us write
(\ref{w3Yw2betagamma}) as
\begin{eqnarray}\label{w3Yw2betagammaexpanded}
\lefteqn{\sum_{l\in{\C}}\sum_{k\in {\mathbb N}}\langle w'_{(3)},
w_{l;\,k}^{\cal Y}w_{(2)}^{(-\beta-\gamma)}\rangle x^{l+1+p}(-1)^k
(\log x)^{k+q}}\nno\\
&&=\sum_{m\in{\C}}\sum_{k\in {\mathbb N}}\langle w'_{(3)},
w_{n-p-1-m;\,k}^{\cal Y}w_{(2)}^{(-\beta-\gamma)}\rangle x^{n-m}(-1)^k
(\log x)^{k+q}
\end{eqnarray}
(recall that we have fixed $n \in \C$).  But each term $\langle
w'_{(3)}, w_{n-p-1-m;\,k}^{\cal Y}w_{(2)}^{(-\beta-\gamma)}\rangle$ in
(\ref{w3Yw2betagammaexpanded}) can be replaced by
\begin{equation}
\langle w'_{(3)}, w_{n-p-1-m;\,k}^{\cal Y}u_{[m]}\rangle,
\end{equation}
where $u_{[m]} \in W_2^{(-\beta-\gamma)}$ is the component of
$w_{(2)}^{(-\beta-\gamma)}$, with respect to the generalized-weight
grading, of (generalized) weight
\[
\wt u_{[m]} = \wt w'_{(3)} - \wt w + n - p - m,
\]
by Proposition \ref{log:logwt}(b).  To see that the coefficient of
$x^n$ in (\ref{w3Yw2betagammaexpanded}) involves only finitely many
powers of $\log x$, independently of the element $w_{(2)}$, we take
$m=0$ in (\ref{w3Yw2betagammaexpanded}) and we observe that the
possible elements $u_{[0]}$ range through the space
\[
(W_2)^{(-\beta -\gamma)}_{[\swt w'_{(3)}-\swt w + n - p]},
\]
which is finite dimensional by the grading restriction condition
(\ref{set:dmfin}).  This proves (\ref{finitelymanypowersoflogx}).  To
prove the lower truncation condition (\ref{log:ltc}), what we must
show is that for sufficiently large $m \in {\mathbb N}$, the
coefficient of $x^{n-m}$ in (\ref{w3Yw2betagammaexpanded}) is $0$
(independently of $w_{(2)}$).  But by the grading restriction
condition (\ref{set:dmltc}),
\[
(W_2)^{(-\beta -\gamma)}_{[\swt w'_{(3)}-\swt w + n - p - m]} = 0
\;\;\mbox{ for }\;m\in {\mathbb N}\;\mbox{ sufficiently large.}
\]
Hence the coefficient of $x^{n-m}$ in (\ref{w3Yw2betagammaexpanded})
is zero for $m\in {\mathbb N}$ sufficiently large, as desired, proving
the lower truncation condition.

For the Jacobi identity, we need to show that
\begin{eqnarray}\label{716}
\lefteqn{\left\langle x^{-1}_0\delta \left( {x_1-x_2\over x_0}\right)
Y_{2}'(v,x_1)A_{r}({\cal
Y})(w_{(1)},x_2)w'_{(3)},w_{(2)}\right\rangle_{W_2}}
\nno \\
&&\;\;\;\;- \left\langle x^{-1}_0\delta \left( {x_2-x_1\over -x_0}\right)
A_{r}({\cal
Y})(w_{(1)},x_2)Y_{3}'(v,x_1)w'_{(3)},w_{(2)}\right\rangle_{W_2}
\nno\\
&&=\left\langle x^{-1}_2\delta \left( {x_1-x_0\over x_2}\right)
A_{r}({\cal Y})(Y_{1}(v,x_0)w_{(1)},x_2)w'_{(3)},w_{(2)}\right
\rangle_{W_2}.
\end{eqnarray}
By the definitions (\ref{y'}) and (\ref{log:Ardef}) we have
\begin{eqnarray}\label{717}
\lefteqn{\langle Y_{2}'(v,x_1)A_{r}({\cal Y})
(w_{(1)},x_2)w'_{(3)},w_{(2)}\rangle_{W_2}}\nno  \\
&&= \langle w'_{(3)},{\cal Y}(e^{x_2L(1)}e^{(2r+1)\pi iL(0)}
(x_2^{-L(0)})^2 w_{(1)},x^{-1}_2)\cdot \nno\\
&&\hspace{8em}\cdot Y_{2}(e^{x_1L(1)}(-x^{-2}_1)^{L(0)}v,x^{-1}_1)w_{(2)}
\rangle_{W_3},
\end{eqnarray}
\begin{eqnarray}\label{718}
\lefteqn{\langle A_{r}({\cal Y})(w_{(1)},x_2)
Y_{3}'(v,x_1)w'_{(3)},w_{(2)}\rangle_{W_2}}\nno  \\
&&=\langle w'_{(3)},Y_{3}(e^{x_1L(1)}(-x^{-2}_1)^{L(0)}v,x^{-1}_1)\cdot \nno\\
&&\hspace{6em}\cdot {\cal Y}(e^{x_2L(1)}e^{(2r+1)\pi iL(0)}
(x_2^{-L(0)})^2 w_{(1)},x^{-1}_2)w_{(2)}\rangle_{W_3},
\end{eqnarray}
\begin{eqnarray}\label{719}
\lefteqn{\langle A_{r}({\cal Y})
(Y_{1}(v,x_0)w_{(1)},x_2)w'_{(3)},w_{(2)}\rangle_{W_2}}\nno \\
&&= \langle w'_{(3)},{\cal Y}(e^{x_2L(1)}e^{(2r+1)\pi iL(0)}
(x_2^{-L(0)})^2
Y_{1}(v,x_0)w_{(1)},x^{-1}_2)w_{(2)}\rangle_{W_3}.
\end{eqnarray}
{From} the Jacobi identity for  ${\cal Y}$  we have
\begin{eqnarray}\label{720}
\lefteqn{\Biggl\langle w'_{(3)},\left( {-x_0\over x_1x_2}\right) ^{-1}\delta
 \left( {x^{-1}_1-x^{-1}_2\over -x_0/x_1x_2} \right)
Y_{3}(e^{x_1L(1)}(-x^{-2}_1)^{L(0)}v,x^{-1}_1)\cdot}\nno \\
&&\;\;\;\;\;\;\;\;\;\;\cdot {\cal Y}(e^{x_2L(1)}e^{(2r+1)\pi iL(0)}
(x_2^{-L(0)})^2 w_{(1)},x^{-1}_2)w_{(2)}\Biggr\rangle_{W_3}\nno \\
&&- \Biggl\langle w'_{(3)},\left( {-x_0\over x_1x_2}\right) ^{-1}\delta
\left( {x^{-1}_2-x^{-1}_1\over x_0/x_1x_2}\right)
{\cal Y}(e^{x_2L(1)}e^{(2r+1)\pi iL(0)}
(x_2^{-L(0)})^2 w_{(1)},x^{-1}_2)\cdot \nno\\
&&\;\;\;\;\;\;\;\;\;\;\cdot
 Y_{2}(e^{x_1L(1)}(-x^{-2}_1)^{L(0)}v,x^{-1}_1)w_{(2)}\Biggr
\rangle_{W_3} \nno\\
&&= \Biggl\langle w'_{(3)},(x^{-1}_2)^{-1}
\delta \left( {x^{-1}_1+x_0/x_1x_2\over x^{-1}_2}\right)
{\cal Y}(Y_{1}(e^{x_1L(1)}(-x^{-2}_1)^{L(0)}v,
-x_0/x_1x_2)\cdot \nno\\
&&\;\;\;\;\;\;\;\;\;\;\cdot e^{x_2L(1)}e^{(2r+1)\pi iL(0)}
(x_2^{-L(0)})^2 w_{(1)},x^{-1}_2)w_{(2)}\Biggr\rangle_{W_3},
\end{eqnarray}
or equivalently,
\begin{eqnarray}\label{721}
\lefteqn{- \Biggl\langle w'_{(3)},x^{-1}_0\delta \left( {x_2-x_1\over -x_0}\right)
Y_{3}(e^{x_1L(1)}(-x^{-2}_1)^{L(0)}v,x^{-1}_1)\cdot }\nno\\
&&\;\;\;\;\;\;\;\;\;\;\cdot {\cal Y}(e^{x_2L(1)}e^{(2r+1)\pi iL(0)}
(x_2^{-L(0)})^2 w_{(1)},
x^{-1}_2)w_{(2)}\Biggr\rangle_{W_3} \nno\\
&&\;\;\;\;+ \Biggl\langle w'_{(3)},x^{-1}_0
\delta \left( {x_1-x_2\over x_0}\right)
{\cal Y}(e^{x_2L(1)}e^{(2r+1)\pi
iL(0)}(x_2^{-L(0)})^2 w_{(1)},x^{-1}_2)\cdot \nno\\
&&\;\;\;\;\;\;\;\;\;\;\cdot Y_{2}(e^{x_1L(1)}(-x^{-2}_1)^{L(0)}v,
x^{-1}_1)w_{(2)}\Biggr\rangle_{W_3} \nno\\
&&= \Biggl\langle w'_{(3)},x^{-1}_1\delta
\left( {x_2+x_0\over x_{1}}\right) {\cal Y}(Y_{1}
(e^{x_1L(1)}(-x^{-2}_1)^{L(0)}v,-x_0/x_1x_2)\cdot\nno\\
&&\;\;\;\;\;\;\;\;\;\;\cdot e^{x_2L(1)}e^{(2r+1)\pi iL(0)}
(x_2^{-L(0)})^2 w_{(1)},x^{-1}_2)w_{(2)}\Biggr\rangle_{W_3}.
\end{eqnarray}
Substituting (\ref{717}), (\ref{718}) and (\ref{719}) into (\ref{716})
and then comparing with (\ref{721}), we see
that the proof of (\ref{716}) is reduced to the proof of the formula
\begin{eqnarray}
\lefteqn{x^{-1}_1\delta \left( {x_2+x_0\over x_1}\right)
{\cal Y}(e^{x_2L(1)}e^{(2r+1)\pi iL(0)}
(x_2^{-L(0)})^2 Y_{1}(v,x_0)w_{(1)},x^{-1}_2)}\nno \\
&&= x^{-1}_1\delta \left( {x_2+x_0\over x_1}\right)
 {\cal Y}(Y_{1}(e^{x_1L(1)}(-x^{-2}_1)^{L(0)}v,-x_0/x_1x_2)\cdot \nno\\
&&\;\;\;\;\;\;\;\;\;\;\cdot e^{x_2L(1)}e^{(2r+1)\pi iL(0)}
(x_2^{-L(0)})^2 w_{(1)},x^{-1}_2),
\end{eqnarray}
or of
\begin{eqnarray}
\lefteqn{{\cal Y}(e^{x_2L(1)}e^{(2r+1)\pi iL(0)}
(x_2^{-L(0)})^2 Y_{1}(v,x_0)w_{(1)},x^{-1}_2)}\nno \\
&&= {\cal Y}(Y_{1}
(e^{(x_2+x_0)L(1)}(-(x_2+x_0)^{-2})^{L(0)}v,-x_0/(x_2+x_0)x_2)\cdot\nno \\
&&\;\;\;\;\;\;\;\;\;\;\cdot e^{x_2L(1)}e^{(2r+1)\pi iL(0)}
(x_2^{-L(0)})^2 w_{(1)},x^{-1}_2).
\end{eqnarray}
We see that we need only prove
\begin{eqnarray}
\lefteqn{e^{x_2L(1)}e^{(2r+1)\pi iL(0)}
(x_2^{-L(0)})^2 Y_{1}(v,x_0)}\nno\\
&&=Y_{1}(e^{(x_2+x_0)L(1)}(-(x_2+x_0)^{-2})^{L(0)}v,
-x_0/(x_2+x_0)x_2)\cdot\nno\\
&&\;\;\;\;\;\;\;\;\;\;\cdot e^{x_2L(1)}e^{(2r+1)\pi iL(0)}
(x_2^{-L(0)})^2
\end{eqnarray}
or equivalently, the conjugation formula
\begin{eqnarray}\label{725}
\lefteqn{e^{x_{2}L(1)}e^{(2r+1)\pi iL(0)}
(x_2^{-L(0)})^2 Y_{1}(v,x_0)(x_2^{L(0)})^2 e^{-(2r+1)\pi iL(0)}
e^{-x_{2}L(1)}}\nno \\
&&= Y_1(e^{(x_{2}+x_0)L(1)}(-(x_{2}+x_0)^{-2})^{L(0)}v,-x_0/(x_{2}+x_0)x_{2})
\end{eqnarray}
for $v\in V$, acting on the module $W_{1}$. But
formula (\ref{725}) follows {from} (\ref{log:p2}), (\ref{log:p3}) and the
formula
\begin{equation}
e^{(2r+1)\pi iL(0)}Y_{1}(v,x)e^{-(2r+1)\pi iL(0)}=Y_{1}((-1)^{L(0)}v, -x),
\end{equation}
which is a special case of (\ref{710}).
This establishes the Jacobi identity.

The $L(-1)$-derivative property follows {}from the same argument used
in the proof of Theorem 5.5.1 of \cite{FHL}: We have (omitting the
subscript $W_3$ on the pairings after a certain point)
\begin{eqnarray}\label{log:ArL(-1)}
\lefteqn{\left\langle \frac{d}{dx}A_r({\cal Y})(w_{(1)},x)w'_{(3)},w_{(2)}
\right\rangle_{W_2}=\frac{d}{dx}\langle A_r({\cal Y})(w_{(1)},x)w'_{(3)},w_{(2)}
\rangle_{W_2}}\nno\\
&&=\frac{d}{dx}\langle w'_{(3)},{\cal Y}(e^{xL(1)}e^{(2r+1)\pi iL(0)}
(x^{-L(0)})^2 w_{(1)},x^{-1})w_{(2)}\rangle_{W_3}\nno\\
&&=\langle w'_{(3)},\frac{d}{dx}{\cal Y}(e^{xL(1)}e^{(2r+1)\pi iL(0)}
(x^{-L(0)})^2 w_{(1)},x^{-1})w_{(2)}\rangle\nno\\
&&=\langle w'_{(3)},{\cal Y}(\frac{d}{dx}(e^{xL(1)}e^{(2r+1)\pi iL(0)}
(x^{-L(0)})^2 w_{(1)}),x^{-1})w_{(2)}\rangle\nno\\
&&\quad+\langle w'_{(3)},\frac{d}{dx}{\cal Y}(w,x^{-1})|_{w=
e^{xL(1)}e^{(2r+1)\pi iL(0)} (x^{-L(0)})^2 w_{(1)}}
w_{(2)}\rangle\nno\\
&&=\langle w'_{(3)},{\cal Y}(e^{xL(1)}L(1)e^{(2r+1)\pi iL(0)}
(x^{-L(0)})^2 w_{(1)},x^{-1})w_{(2)}\rangle\nno\\
&&\quad+\langle w'_{(3)},{\cal Y}(e^{xL(1)}(-2L(0)x^{-1})
e^{(2r+1)\pi iL(0)}(x^{-L(0)})^2w_{(1)},x^{-1})w_{(2)}\rangle\nno\\
&&\quad+\langle w'_{(3)},\frac{d}{d{x^{-1}}}{\cal Y}(w,x^{-1})|_{w=
e^{xL(1)}e^{(2r+1)\pi iL(0)} (x^{-L(0)})^2 w_{(1)}}
w_{(2)}\rangle(-x^{-2})\nno\\
&&=\langle w'_{(3)},{\cal Y}(e^{xL(1)}(2xL(0)-x^2L(1))e^{(2r+1)\pi iL(0)}
(x^{-L(0)})^2(-x^{-2})w_{(1)},x^{-1})w_{(2)}\rangle\nno\\
&&\quad+\langle w'_{(3)},{\cal Y}(L(-1)e^{xL(1)}e^{(2r+1)\pi iL(0)}
(x^{-L(0)})^2 w_{(1)},x^{-1})w_{(2)}\rangle(-x^{-2})\nno\\
&&=\langle w'_{(3)},{\cal Y}(e^{xL(1)}(2xL(0)-x^2L(1))e^{(2r+1)\pi iL(0)}
(x^{-L(0)})^2(-x^{-2})w_{(1)},x^{-1})w_{(2)}\rangle\nno\\
&&\quad+\langle w'_{(3)},{\cal Y}(L(-1)e^{xL(1)}e^{(2r+1)\pi iL(0)}
(x^{-L(0)})^2 (-x^{-2})w_{(1)},x^{-1})w_{(2)}\rangle
\end{eqnarray}
Now by (\ref{log:SL2-3}), with $x$ replaced by $-x$, we have
\[
L(-1)e^{xL(1)}=e^{xL(1)}(L(-1)-2xL(0)+x^2L(1)).
\]
Using this together with (\ref{log:xLx^}) and (\ref{log:SL2-2}) (with
$x$ specialized to $-(2r+1)\pi i$, and the convergence of the
exponential series invoked; recall Remark
\ref{analyticallyconvergent}), we see that the right-hand side of
(\ref{log:ArL(-1)}) equals
\begin{eqnarray*}
\lefteqn{\langle w'_{(3)},{\cal Y}(e^{xL(1)}L(-1)e^{(2r+1)\pi iL(0)}
(x^{-L(0)})^2(-x^{-2})w_{(1)},x^{-1})w_{(2)}\rangle}\\
&&=\langle w'_{(3)},{\cal Y}(e^{xL(1)}e^{(2r+1)\pi iL(0)}(-L(-1))
(x^{-L(0)})^2(-x^{-2})w_{(1)},x^{-1})w_{(2)}\rangle\\
&&=\langle w'_{(3)},{\cal Y}(e^{xL(1)}e^{(2r+1)\pi iL(0)}
(x^{-L(0)})^2L(-1)w_{(1)},x^{-1})w_{(2)}\rangle\\
&&=\langle A_r({\cal Y})(L(-1)w_{(1)},x)w'_{(3)},w_{(2)}\rangle_{W_2},
\end{eqnarray*}
as desired.

We now show that, in case $V$ is M\"obius, the ${\mathfrak
s}{\mathfrak l}(2)$-bracket relations (\ref{log:L(j)b}) hold for
$A_r({\cal Y})$.  For these, we first see that, for $j=-1,0,1$, by
using (\ref{L'(n)}), (\ref{log:Ardef}) and the ${\mathfrak
s}{\mathfrak l}(2)$-bracket relations (\ref{log:L(j)b}) for ${\cal Y}$
we have
\begin{eqnarray}\label{log:tmp1}
\lefteqn{\langle L'(j)A_r({\cal Y})(w_{(1)},x)w'_{(3)},w_{(2)}
\rangle_{W_2}}\nno\\
&&=\langle A_r({\cal Y})(w_{(1)},x)w'_{(3)},L(-j)w_{(2)}
\rangle_{W_2}\nno\\
&&=\langle w'_{(3)},{\cal Y}(e^{xL(1)}e^{(2r+1)\pi iL(0)}(x^{-L(0)})^2
w_{(1)},x^{-1})L(-j)w_{(2)}\rangle_{W_3}\nno\\
&&=\langle w'_{(3)},L(-j){\cal Y}(e^{xL(1)}e^{(2r+1)\pi iL(0)}
(x^{-L(0)})^2w_{(1)},x^{-1})w_{(2)}\rangle_{W_3}\nno\\
&&\hspace{1em}-\Biggl\langle w'_{(3)},\sum_{i=0}^{-j+1}{-j+1\choose i}x^{-i}
\cdot\nno\\
&&\hspace{8em}\cdot{\cal Y}(L(-j-i)e^{xL(1)}e^{(2r+1)\pi iL(0)}
(x^{-L(0)})^2w_{(1)},x^{-1})w_{(2)}\Biggr\rangle_{W_3}.
\end{eqnarray}
Now {}from (\ref{log:SL2-3}), (\ref{log:SL2-2}) (with $x$ specialized to
$-(2r+1)\pi i$) and (\ref{log:xLx^}), one computes that
\begin{eqnarray*}
\lefteqn{(x^{L(0)})^2e^{-(2r+1)\pi iL(0)}e^{-xL(1)}
\left(\begin{array}{c}L(-1)\\L(0)\\L(1)\end{array}\right)
e^{xL(1)}e^{(2r+1)\pi iL(0)}(x^{-L(0)})^2}\\
&&\hspace{5cm}=\left(\begin{array}{ccc}-x^2&-2x&-1\\0&1&x^{-1}
\\0&0&-x^{-2}\end{array}\right)
\left(\begin{array}{c}L(-1)\\L(0)\\L(1)\end{array}\right)
\end{eqnarray*}
on $W_1$, which implies that
\begin{eqnarray*}
\lefteqn{\sum_{i=0}^{-j+1}{-j+1\choose i}x^{-i}L(-j-i)e^{xL(1)}
e^{(2r+1)\pi iL(0)}(x^{-L(0)})^2}\\
&&\hspace{2em}=-\sum_{i=0}^{j+1}{j+1\choose i}x^ie^{xL(1)}e^{(2r+1)
\pi iL(0)}(x^{-L(0)})^2L(j-i)
\end{eqnarray*}
on $W_1$.  Hence the right-hand side of
(\ref{log:tmp1}) is equal to
\begin{eqnarray*}
&&\lefteqn{\langle L'(j)w'_{(3)},{\cal Y}(e^{xL(1)}e^{(2r+1)\pi iL(0)}
(x^{-L(0)})^2w_{(1)},x^{-1})w_{(2)}\rangle_{W_3}}\\
&&\hspace{1em}+\sum_{i=0}^{j+1}{j+1\choose i}x^i\langle w'_{(3)},
{\cal Y}(e^{xL(1)}e^{(2r+1)\pi iL(0)}(x^{-L(0)})^2
L(j-i)w_{(1)},x^{-1})w_{(2)} \rangle_{W_2}\\
&&=\langle A_r({\cal Y})(w_{(1)},x)L'(j)w'_{(3)},w_{(2)}
\rangle_{W_2}\\
&&\hspace{1em}+\sum_{i=0}^{j+1}{j+1\choose i}x^i\langle A_r({\cal Y})
(L(j-i)w_{(1)},x)w'_{(3)},w_{(2)} \rangle_{W_2},
\end{eqnarray*}
and the ${\mathfrak s}{\mathfrak l}(2)$-bracket relations for
$A_r({\cal Y})$ are proved.  (Note that this argument essentially
generalizes the proof of Lemma \ref{sl2opposite}.)  We have finished
proving that $A_r({\cal Y})$ is a grading-compatible logarithmic
intertwining operator.

Finally, for the relation (\ref{log:ar}), we of course identify
$W''_2$ with $W_2$ and $W''_3$ with $W_3$, according to Theorem
\ref{set:W'}.  Let us view ${\cal Y}$ as a grading-compatible
logarithmic intertwining operator of type ${W'_2\choose W_1\,W'_3}$,
so that $A_r({\cal Y})$ is such an operator of type ${W_3\choose
W_1\,W_2}$.  We have
\begin{eqnarray*}
\lefteqn{\langle A_{-r-1}A_r({\cal Y})(w_{(1)},x)w'_{(3)},w_{(2)}
\rangle_{W_2}}\\
&&=\langle w'_{(3)},A_r({\cal Y})(e^{xL(1)}e^{(-2r-1)\pi iL(0)}
(x^{-L(0)})^2 w_{(1)},x^{-1})w_{(2)}\rangle_{W_3}\\
&&=\langle {\cal Y}(e^{x^{-1}L(1)}e^{(2r+1)\pi iL(0)}(x^{L(0)})^2
e^{xL(1)}e^{(-2r-1)\pi iL(0)}(x^{-L(0)})^2w_{(1)},x)w'_{(3)},w_{(2)}
\rangle_{W_2}\\
&&=\langle {\cal Y}(w_{(1)},x)w'_{(3)},w_{(2)}\rangle_{W_2},
\end{eqnarray*}
where the last equality is due to the relation
\begin{equation}\label{conjrelation}
e^{(2r+1)\pi iL(0)}(x^{L(0)})^2 e^{xL(1)}(x^{-L(0)})^2 e^{-(2r+1)\pi
iL(0)} = e^{-x^{-1}L(1)}
\end{equation}
on $W_{1}$, whose proof is similar to that of formula (5.3.1) of
\cite{FHL}.  Namely, (\ref{conjrelation}) follows {}from the relation
\begin{equation}\label{xto-1/x}
e^{(2r+1)\pi iL(0)}(x^{L(0)})^2 xL(1) (x^{-L(0)})^2 e^{-(2r+1)\pi
iL(0)} = -x^{-1}L(1),
\end{equation}
which realizes the transformation $x \mapsto -\frac{1}{x}$, and 
(\ref{xto-1/x}) follows {}from (\ref{log:xLx^}) together with 
(\ref{log:SL2-2}) specialized as above.  \epf

\begin{rema}{\rm
The last argument in the proof shows that for any $r,s\in{\mathbb Z}$,
formula (\ref{log:ar}) generalizes to:
\[
A_s(A_r({\cal Y}))={\cal Y}_{[0,r+s+1,0]}
\]
(recall Remark \ref{Ys1s2s3}).
}
\end{rema}

With $V$, $W_1$, $W_2$ and $W_3$ strongly graded, set
\begin{equation}
N_{W_1W_2W_3}=N_{W_1\,W_2}^{W'_3}.
\end{equation}
Then Proposition \ref{log:omega} gives 
\[
N_{W_1W_2W_3}=N_{W_2W_1W_3}
\]
and Proposition \ref{log:A} gives
\[
N_{W_1W_2W_3}=N_{W_1W_3W_2}.
\]
Thus for any permutation $\sigma$ of $(1,2,3)$,
\begin{equation}
N_{W_1W_2W_3}=N_{W_{\sigma(1)}W_{\sigma(2)}W_{\sigma(3)}}.
\end{equation}

It is clear {}from Proposition \ref{log:logwt}(b) that in the nontrivial
logarithmic intertwining operator case, taking projections of ${\cal
Y}(w_{(1)},x)w_{(2)}$ to (generalized) weight subspaces is not enough
to recover its coefficients of $x^n(\log x)^k$ for $n\in{\mathbb C}$
and $k\in{\mathbb N}$, in contrast with the (ordinary) intertwining
operator case (cf.\ \cite{tensor1}, the paragraph containing formula
(4.17)).  However, taking projections of certain related intertwining
operators does indeed suffice for this purpose:

\begin{propo}\label{log:proj}
Let $W_1$, $W_2$, $W_3$ be generalized modules for a M\"obius (or
conformal) vertex algebra $V$ and let ${\cal Y}$ be a logarithmic
intertwining operator of type ${W_3\choose W_1\,W_2}$.  Let
$w_{(1)}\in W_1$ and $w_{(2)}\in W_2$ be homogeneous of generalized
weights $n_1$ and $n_2$, respectively. Then for any $n\in {\mathbb C}$
and any $r\in {\mathbb N}$, ${w_{(1)}}^{\cal Y}_{n;\,r}w_{(2)}$ can be
written as a certain linear combination of products of the component of 
weight $n_1+n_2-n-1$ of
\begin{equation}\label{modified-lio}
(L(0)-n_1-n_2+n+1)^l{\cal Y}((L(0)-n_1)^iw_{(1)},x)(L(0)-n_2)^jw_{(2)}
\end{equation}
for certain $i,j,l\in{\mathbb N}$ with monomials of the form
$x^{n+1}(\log x)^m$ for certain $m\in{\mathbb N}$.
\end{propo}
\pf Multiplying (\ref{log:r+t=?}) by $x^{-n-1}(\log x)^k$ and summing
over $k\in {\mathbb N}$ (a finite sum by definition) we have that for
any $t\in {\mathbb N}$,
\begin{eqnarray}\label{log:proj1}
\lefteqn{\sum_{k\in {\mathbb N}}{k+t\choose t}{w_{(1)}}^{\cal
Y}_{n;\,k+t}w_{(2)}x^{-n-1}(\log x)^k \;
\bigg(=x^{-n-1}\sum_{k\in {\mathbb N}}{k\choose t}
{w_{(1)}}^{\cal Y}_{n;\,k}w_{(2)}(\log x)^{k-t}\bigg)}\nno\\
&&=\sum_{i,j,l\in {\mathbb N}, \, i+j+l=t}\frac{1}{i!j!l!}(-1)^{i+j}
\sum_{k\in {\mathbb N}}(L(0)-n_1-n_2+n+1)^l\cdot\nno\\
&&\quad\cdot((L(0)-n_1)^i{w_{(1)}})^{\cal Y}_{n;\,k}
(L(0)-n_2)^jw_{(2)})x^{-n-1}(\log x)^k\nno\\
&&=\sum_{i,j,l\in {\mathbb N}, \, i+j+l=t}\frac{1}{i!j!l!}(-1)^{i+j}
\pi_{n_1+n_2-n-1}((L(0)-n_1-n_2+n+1)^l\cdot\nno\\
&&\quad\cdot{\cal Y}((L(0)-n_1)^i{w_{(1)}},x) (L(0)-n_2)^jw_{(2)}).
\end{eqnarray}
Let $K$ be a positive integer such that ${w_{(1)}}^{\cal Y}_{n;\,k'}
w_{(2)}=0$ for all $k'\geq K$, Denote the right-hand side of
(\ref{log:proj1}) by $\pi(t,w_{(1)},w_{(2)},x,\log x)$.  Then by
putting the identities (\ref{log:proj1}) for $t=0,1,\dots,K-1$ together in
matrix form we have
\begin{equation}\label{log:projmat}
x^{-n-1}A\left(\begin{array}{c}{w_{(1)}}^{\cal Y}_{n;\,0}w_{(2)}\\
{w_{(1)}}^{\cal Y}_{n;\,1}w_{(2)}\\\vdots\\
{w_{(1)}}^{\cal Y}_{n;\,K-1}w_{(2)}\end{array}\right)=
\left(\begin{array}{c}\pi(0,w_{(1)},w_{(2)},x,\log x)\\
\pi(1,w_{(1)},w_{(2)},x,\log x)\\\vdots\\
\pi(K-1,w_{(1)},w_{(2)},x,\log x)\end{array}\right)
\end{equation}
where $A$ is the $K\times K$ matrix whose $(i,j)$-entry is equal to
$\displaystyle{j-1\choose i-1}(\log x)^{j-i}$.  Letting $P_K$ be the
triangular matrix whose $(i,j)$-entry is $\displaystyle{j-1\choose
i-1}$ (an upper triangular ``Pascal matrix''), we have
\[
A={\rm diag}(1,(\log x)^{-1},\dots,(\log x)^{-(K-1)})\cdot P_K\cdot
{\rm diag}(1,\log x,\dots,(\log x)^{K-1}).
\]
Its inverse is
\[
A^{-1}={\rm diag}(1,(\log x)^{-1},\dots,(\log x)^{-(K-1)})\cdot
P_K^{-1}\cdot{\rm diag}(1,\log x,\dots,(\log x)^{K-1})
\]
and the $(i,j)$-entry of $P_K^{-1}$ is $\displaystyle(-1)^{i+j}
{j-1\choose i-1}$.  Now multiplying the left-hand side of
(\ref{log:projmat}) by $x^{n+1}A^{-1}$ we obtain
\[
\left(\begin{array}{c}{w_{(1)}}^{\cal Y}_{n;\,0}w_{(2)}\\
{w_{(1)}}^{\cal Y}_{n;\,1}w_{(2)}\\\vdots\\
{w_{(1)}}^{\cal Y}_{n;\,K-1}w_{(2)}\end{array}\right)=x^{n+1}A^{-1}
\left(\begin{array}{c}\pi(0,w_{(1)},w_{(2)},x,\log x)\\
\pi(1,w_{(1)},w_{(2)},x,\log x)\\\vdots\\
\pi(K-1,w_{(1)},w_{(2)},x,\log x)\end{array}\right)
\]
or explicitly,
\begin{equation}\label{log:last}
(w_{(1)})^{\cal Y}_{n;\,r}w_{(2)}=x^{n+1}\sum_{t=r}^{K-1}(-1)^{r+t}
{t\choose r}(\log x)^{t-r}\pi(t,w_{(1)},w_{(2)},x,\log x)
\end{equation}
for $r=0,1,\dots,K-1$. (In particular, all $x$'s and $\log x$'s cancel
out in the right-hand side of (\ref{log:last}).)  \epfv

\begin{rema}{\rm Proposition \ref{log:proj} says that
taking projections of (\ref{modified-lio}) to (generalized) weight
subspaces is enough to recover the coefficients of $\Y$.  Although we
will not need to use Proposition \ref{log:proj} in the proof of our
main results, we will heavily use certain analytic results, in
particular, Propositions \ref{real-exp-set} and
\ref{log-coeff-conv<=>iterate-conv} and Corollary
\ref{double-conv<=>iterate-conv} in Section 7, to determine, using
projections, coefficients of products of powers of complex variables
and of their logs, in expressions involving logarithmic intertwining
operators, in many results that we will need, including in particular
Propositions \ref{prod=0=>comp=0} and \ref{iter=0=>comp=0} and their
respective corollaries.}
\end{rema}


\bigskip

\noindent {\small \sc Department of Mathematics, Rutgers University,
Piscataway, NJ 08854 (permanent address)}

\noindent {\it and}

\noindent {\small \sc Beijing International Center for Mathematical Research,
Peking University, Beijing, China}

\noindent {\em E-mail address}: yzhuang@math.rutgers.edu

\vspace{1em}

\noindent {\small \sc Department of Mathematics, Rutgers University,
Piscataway, NJ 08854}

\noindent {\em E-mail address}: lepowsky@math.rutgers.edu

\vspace{1em}

\noindent {\small \sc Department of Mathematics, Rutgers University,
Piscataway, NJ 08854}

\noindent {\em E-mail address}: linzhang@math.rutgers.edu

\end{document}